\documentclass{amsart}
\usepackage{amsmath,amssymb} 
\usepackage{amscd}
\usepackage[dvipdfm]{graphicx}

\theoremstyle{plain}
\newtheorem{thm}{Theorem}[section]
\newtheorem{lem}{Lemma}[section]

\newtheorem{prop}{Proposition}[section]
\newtheorem{clm}{Claim}[section]
\theoremstyle{definition}
\newtheorem{defn}{Definition}[section]
\theoremstyle{remark}
\newtheorem*{rem}{Remark}
\newtheorem*{prf}{Proof}

\newtheorem*{prf1}{Proof of Lemma \ref{lem1}}
\newtheorem*{prf2}{Proof of Lemma \ref{lem2}}
\newtheorem*{prf3}{Proof of Lemma \ref{lem3}}
\newtheorem*{prf4}{Proof of Proposition \ref{proplast}}
\newtheorem*{prf5}{Proof of Theorem \ref{thmMain}}



\title{Heegaard Floer homology, L-spaces and smoothing order on links I}
\author{TAKUYA USUI}

\begin{document}

\begin{abstract}
In this paper, we focus on L-spaces for which the boundary maps of the Heegaard Floer chain complexes vanish. We collect such manifolds systematically by using the smoothing order on links. 
\end{abstract}

\address{Graduate School of Mathematical Science, University of Tokyo, 3-8-1 Komaba Meguroku Tokyo 153-8914, Japan}
\email{t.usuiusu.t@gmail.com}

\maketitle

%
\section{Introduction}

\ In \cite{OS2} and \cite{OS1}, Ozsv\'{a}th and Szab\'{o} introduced the \textit{Heegaard-Floer homology} $\widehat{HF}(Y)$ for a closed oriented three manifold $Y$. The Heegaard Floer homology $\widehat{HF}(Y)$ is defined by using a pointed Heegaard diagram representing $Y$ and a certain version of Lagrangian Floer theory. The boundary map of the chain complex \textit{counts} the number of pseudo-holomorphic Whitney disks. Of course, the boundary map depends on the pointed Heegaard diagram. In this paper, the coefficient of homology is $\mathbb{Z}_{2}$. A rational homology three-sphere $Y$ is called an L-space when its Heegaard Floer homology $\widehat{HF}(Y)$ is a $\mathbb{Z}_{2}$-vector space with dimension $|H_{1}(Y;\mathbb{Z})|$, where $|H_{1}(Y;\mathbb{Z})|$ is the number of elements in $H_{1}(Y;\mathbb{Z})$. 


In this paper, we consider a special class of L-spaces. 
\begin{defn}
An L-space $Y$ is \textit{strong} if there is a pointed Heegaard diagram representing $Y$ such that the boundary map vanishes.
\end{defn}

Strong L-spaces are originally defined in \cite{LL} in another way (see Proposition \ref{equi}), and discussed in \cite{BGW} and \cite{Greene}.
 
Now, We prepare some notations to state the main theorems.

\begin{figure}[h]

\includegraphics[width=7cm,clip]{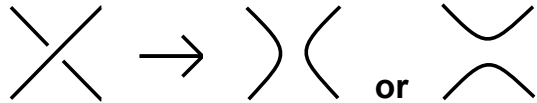} \\
\caption{smoothing}\label{smoothing}
\end{figure}

For a link $L$ in $S^3$, we can get a link diagram $D_{L}$ in $S^2$ by projecting $L$ to $S^2 \subset S^3$. To make other link diagrams from $D_{L}$, we can smooth a crossing point in different two ways (see Figure \ref{smoothing}.)

In \cite{EIT} and \cite{Taniyama}, the following ordering on links is defined.
\begin{defn}

Let $D_{L_{1}}$ and $D_{L_{2}}$ be alternating link diagrams in $S^2$. 
We say $D_{L_{1}} \subseteq D_{L_{2}}$ if $D_{L_{2}}$ contains $D_{L_{1}}$ as a connected component after smoothing some crossing points of $D_{L_{2}}$.

Let $L_{1}$ and $L_{2}$ be alternating links in $S^3$. Then, we say $L_{1} \leq L_{2}$ if for any minimal crossing alternating link diagram $D_{L_{2}}$ of $L_{2}$, there is a minimal crossing alternating link diagram $D_{L_{1}}$ of $L_{1}$ such that $D_{L_{1}} \subseteq D_{L_{2}}$. 
\end{defn}

\begin{figure}[h]

\includegraphics[width=4cm,clip]{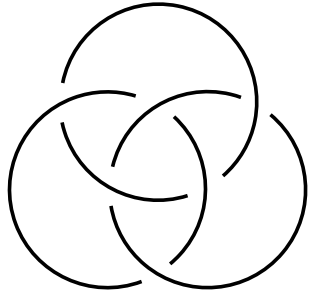} \\
\caption{The Borromean rings}\label{borromean}
\end{figure}

These orderings on links and diagrams are called \textit{smoothing orders} in \cite{EIT}. Note that smoothing orders become partial orderings. Let us denote the minimal crossing number of $L$ by $c(L)$. If $L_{1} \leq L_{2}$, then $c(L_{1}) \leq c(L_{2})$. We can check the well-definedness by using this observation. Actually, if $L_{1} \leq L_{2}$ and $L_{2} \leq L_{1}$, then $c(L_{1}) = c(L_{2})$ and there is no smoothed crossing point. So $L_{1} = L_{2}$. Next, if $L_{1} \leq L_{2}$ and $L_{2} \leq L_{3}$, then $L_{1} \leq L_{3}$ by defintion.
Note that we can define $\leq$ for any two links by ignoring alternating conditions. But in this paper we consider only alternating links and alternating link diagrams. The Borromean rings $\rm{Brm}$ are an alternating link in $S^3$ whose diagram looks  as in Figure \ref{borromean}. We fix this diagram and denote it by $\rm{Brm}$ too.

\begin{defn}
$\mathcal{L}_{\overline{\rm{Brm}}} = \{$ an alternating link $L$ in $S^3$ such that $\rm{Brm} \nleq L\}$, where $\rm{Brm}$ is the Borromean rings.
\end{defn}

Denote $\Sigma(L)$ a double branched covering of $S^3$ branched along a link $L$.
The first main result is as follows:

\begin{thm} \label{thmMain}
Let $L$ be a link in $S^3$. If $L$ satisfies the following conditions:
\begin{itemize}
\item $L \in \mathcal{L}_{\overline{\rm{Brm}}}$, 
\item $\Sigma(L)$ is a rational homology three-sphere,
\end{itemize}
then $\Sigma(L)$ is a strong L-space and a graph manifold (or a connected sum of graphmanifolds).
\end{thm}

A graph manifold is defined as follows.
\begin{defn}
A closed oriented three manifold $Y$ is a \textit{graph manifold} if $Y$ can be decomposed along embedded tori into finitely many Seifert manifolds.
\end{defn}

Now, we recall the following fact. It is proved in \cite{Greene}.
\begin{thm}\label{strongthm}
For an alternating link $L$, if $\Sigma(L)$ is a raional homology three-sphere, $\Sigma(L)$ becomes a strong L-space.
\end{thm}

Theorem \ref{strongthm} seems stronger than Theorem \ref{thmMain}. But, we prove Theorem \ref{thmMain} in a different way. Moreover, we collect systematically strong L-spaces which become graph manifolds.





\section{Heegaard-Floer homology and L-spaces}


\ The Heegaard Floer homology of a closed oriented three manifold $Y$ is defined from a pointed Heegaard diagram representing $Y$.
Let $f$ be a self-indexing Morse function on $Y$ with $1$ index zero critical point and $1$ index three critical point. Then, $f$ gives a Heegaard splitting of $Y$. That is, $Y$ is given by glueing two handlebodies $f^{-1}([0,3/2])$ and $f^{-1}([3/2, 3])$ along their boundaries. If the number of index one critical points or the number of index two critical points of $f$ is $g$, then $\Sigma = f^{-1}(3/2)$ is a closed oriented genus $g$ surface. 
We fix a gradient flow on $Y$ corresponding to $f$. We get a collection $\alpha = \{ \alpha_{1}, \cdots, \alpha_{g} \}$ of $\alpha$ curves on $\Sigma$ which flow down to the index one critical points, and another collection $\beta = \{ \beta_{1}, \cdots, \beta_{g} \}$ of $\beta$ curves on $\Sigma$ which flow up to the index two critical points.
Let $z$ be a point in $\Sigma \setminus (\alpha \cup \beta)$.
The tuple $(\Sigma,\alpha,\beta,z)$ is called a \textit{pointed Heegaard diagram} for $Y$.
Note that $\alpha$ and $\beta$ curves are characterized as pairwise disjoint, homologically linearly independent, simple closed curves on $\Sigma$. We can assume $\alpha$-curves intersect $\beta$-curves transversaly.

Next, we review the definition of the Heegaard Floer chain complex.

Let $(\Sigma,\alpha,\beta,z)$ be a pointed Heegaard diagram for $Y$.
The $g$-fold \textit{symmetric product} of the closed oriented surface $\Sigma$ is defined by $\rm{Sym}^{g}(\Sigma) = \Sigma^{\times g}/S_{g}$. That is, the quotient of $\Sigma^{\times g}$ by the natural action of the symmetric group on $g$ letters. 





Let us define $\mathbb{T}_{\alpha } = \alpha_{1} \times \cdots \times \alpha_{g} /S_{g}$ and $\mathbb{T}_{\beta} = \beta_{1} \times \cdots \times \beta_{g}/S_{g}$.
Then, the chain complex $\widehat{CF}(\Sigma,\alpha,\beta,z)$ is defined as a $\mathbb{Z}_{2}$-vector space generated by the elements of 
$$\mathbb{T}_{\alpha } \cap \mathbb{T}_{\beta }= \{ x=(x_{1 \sigma (1)},x_{2 \sigma (2)},\cdots , x_{g \sigma (g)}) | x_{i \sigma (i)} \in \alpha_{i} \cap \beta_{\sigma(i)}, \sigma \in S_{g} \}.$$ 

Then, the boundary map $\widehat{\partial}$ is given by

\begin{equation}\label{def}
\widehat{\partial}x = \sum_{y \in \mathbb{T}_{\alpha} \cap \mathbb{T}_{\beta}} c(x,y) \cdot y,
\end{equation}
where $c(x,y) \in \mathbb{Z}_{2}$ is defined by \textit{counting} the number of pseudo-holomorphic Whitney disks. For more details, see \cite{OS2}. 






\begin{defn}\cite{OS2}
The homology of the chain complex $(\widehat{CF}(\Sigma,\alpha,\beta,z),\widehat{\partial})$ is called
the Heegaard Floer homology of a pointed Heegaard diagram. We denote it by $\widehat{HF}(\Sigma,\alpha,\beta,z)$.
\end{defn}

\begin{rem}

For appropriate pointed Heegaard diagrams representing $Y$, their Heegaard Floer homologies become isomorphic. So we can define the Heegaard Floer homology of $Y$. Denote it by $\widehat{HF}(Y)$. (For more details, see \cite{OS2}).

\end{rem}


In this paper, we consider only L-spaces, in particular strong L-spaces. The following proposition enables us to define strong L-spaces in another way. The second condition comes from \cite{LL}.

\begin{prop}\label{equi}
Let $(\Sigma,\alpha,\beta,z)$ be a pointed Heegaard diagram representing a rational homology sphere $Y$. Then, the following two conditions (1) and (2) are equivalent. 
\begin{enumerate}
\item the boundary map $\widehat{\partial}$ is the zero map, and $Y$ is an $L$-apace.
\item $|\mathbb{T}_{\alpha} \cap \mathbb{T}_{\beta}|=|H_{1}(Y;\mathbb{Z})|$.
\end{enumerate}
\end{prop}

For example, any lens-spaces are strong L-spaces. Actually, we can draw a genus one Heegaard diagram representing $L(p,q)$ for which the two circles $\alpha$ and $\beta$ meet transversely in $p$ points. That is, $|\mathbb{T}_{\alpha} \cap \mathbb{T}_{\beta}|=|H_{1}(L(p,q);\mathbb{Z})| = p$.


To prove this proposition, we recall that the Heegaard Floer homology $\widehat{HF}(Y)$ admits a relative $\mathbb{Z}/ 2 \mathbb{Z}$ grading(\cite{OS1}) By using this grading, the Euler characteristic satisfies the following equation.
$$\chi(\widehat{HF}(Y,s)) = |H_{1}(Y;\mathbb{Z})|.$$




\begin{prf}
The first condition tells us that $\widehat{CF}(Y)$ becomes a $\mathbb{Z}_{2}$-vector space with dimension $|H_{1}(Y;\mathbb{Z})|$. By definition of $\widehat{CF}(Y)$, we get that $|\mathbb{T}_{\alpha} \cap \mathbb{T}_{\beta}|=|H_{1}(Y;\mathbb{Z})|$.
Conversely, the second condition and the above equation tell us that both $\widehat{CF}(Y)$ and $\widehat{HF}(Y)$
become $\mathbb{Z}_{2}$-vector spaces with dimension $|H_{1}(Y;\mathbb{Z})|$. Therefore, the first condition follows.
\qed
\end{prf}

%
\section{$B$-reducible alternating links and Smoothing order}

\ In this section, we introduce some link type specializing alternating links by using the smoothing order. Arfer that, we prove that the link type is the same as $\mathcal{L}_{\overline{\rm{Brm}}}$.

\subsection{$B$-reducible alternating links}

Let us denote $D_{alt}$ the set of alternating link diagrams in $\mathbb{R}^2$ modulo isotopies.

\begin{figure}[h]

\includegraphics[width=6cm,clip]{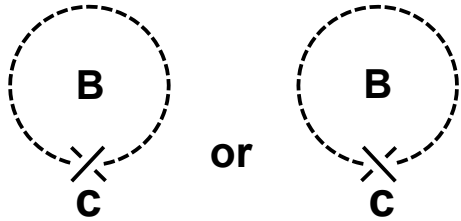} \\
\caption{1-reducible}\label{1-reducible}
\end{figure}

\begin{figure}[h]

\includegraphics[width=6cm,clip]{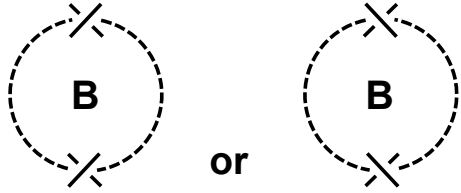} \\
\caption{2-reducible}\label{2-reducible}
\end{figure}

\begin{defn}
Let $D_{L}$ be in $D_{alt}$. An embedded disk $B$ in $\mathbb{R}^2$ is called 1-reducible for $D_{L}$ if the boundary of $B$ intersects with $D_{L}$ at just one crossing point $c$ and $c$ looks  as in Figure \ref{1-reducible}. Similarly $B$ is called 2-reducible for $D_{L}$ if the boundary of $B$ intersects with $D_{L}$ at just two crossing points $c_{1}$ and $c_{2}$ and they look as in Figure \ref{2-reducible}.
In short, $B$ is called reducible for $D_{L}$ if it is 1- or 2-reducible for $D_{L}$.
\end{defn}

For a reducible disk $B$ for $D_{L}$, we define some operations and then get new alternating link diagrams as follows.

\begin{figure}[h]

\includegraphics[width=6cm,clip]{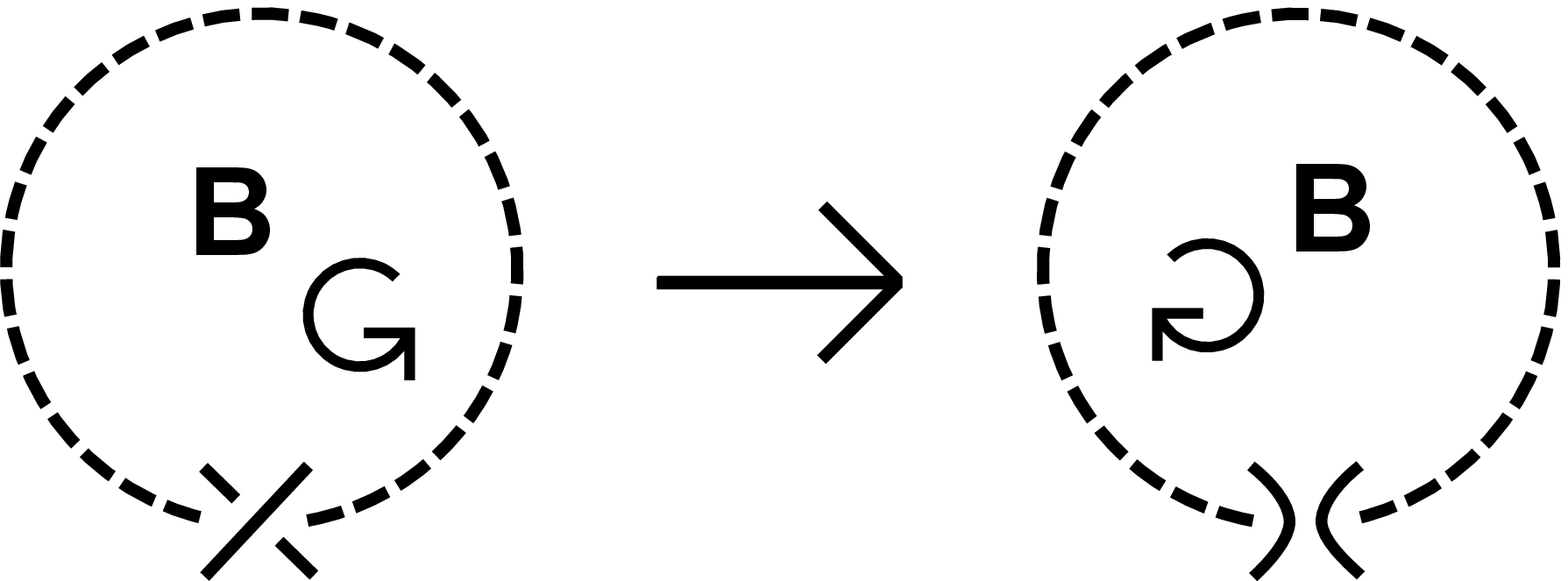} \\
\caption{(I)-move}\label{move1}
\end{figure}

\begin{figure}[h]

\includegraphics[width=6cm,clip]{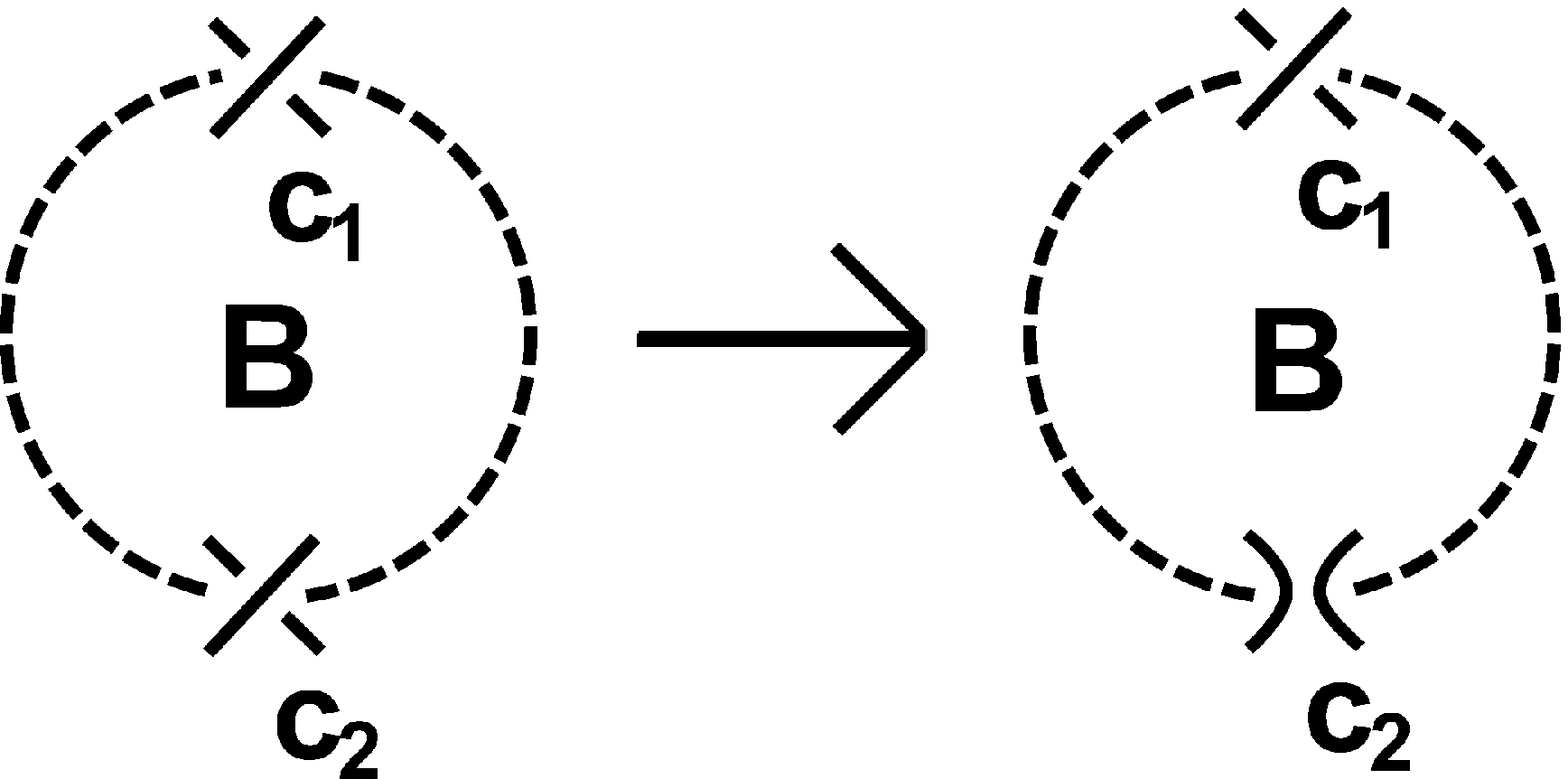} \\
\caption{(II)-move}\label{move2}
\end{figure}

\begin{defn}

If there is a 1-reducible disk $B$ for $D_{L}$, we can get a new alternating diagram $D_{L}(B)$ by reversing the disk $B$ together with the link in $B$ to eliminate the crossing point $c$ (see Figure \ref{move1}). This is called (I)-move.
If there is a 2-reducible disk $B$ for $D_{L}$, we can get two possible alternating link diagrams by smoothing one of the two crossing points $c_{1}$ and $c_{2}$ as in Figure \ref{move2}. We call these two diagrams $D_{L}(B)$ without distinction. This is called (II)-move.
In short, we can get a new alternating link diagram $D_{L}(B)$ in $D_{alt}$ from $D_{L}$ and $B$ by using one of the above operations.

\end{defn}

Now we define a subclass of alternating link diagrams.

\begin{defn}
A class $D_{red}$ is defined as the subset of $D_{alt}$ whose element $D_{L}$ satisfies the one of the following two properties.
\begin{itemize}
\item $D_{L}$ is a disjoint union of finite number of the unknot diagrams. 
\item $D_{L}$ is not a disjoint union of finite number of the unknot diagrams, but there are a sequence of embedded disks ${B_{1},\cdots ,B_{n}}$ and a sequence of (I) or (II)-moves such that 
 \begin{itemize}
 \item $B_{1}$ is reducible for $D_{L}$, 
 \item $B_{2}$ is reducible for $D_{L}(B_{1})$, 
 \item $B_{3}$ is reducible for $D_{L}(B_{1},B_{2}) = D_{L}(B_{1})(B_{2})$,

$\vdots$

 \item $B_{n}$ is reducible for $D_{L}(B_{1},\cdots ,B_{n-1})$,
 \item $D_{L}(B_{1},\cdots ,B_{n})$ is a disjoint union of finite number of the unknot diagrams.
 \end{itemize}
 Note that the expressions $D_{L}(B_{1},\cdots ,B_{m})$ depend on the choice of the operations if the reducible disks are 2-reducible.
\end{itemize}
\end{defn}

\begin{figure}[h]

\includegraphics[width=8cm,clip]{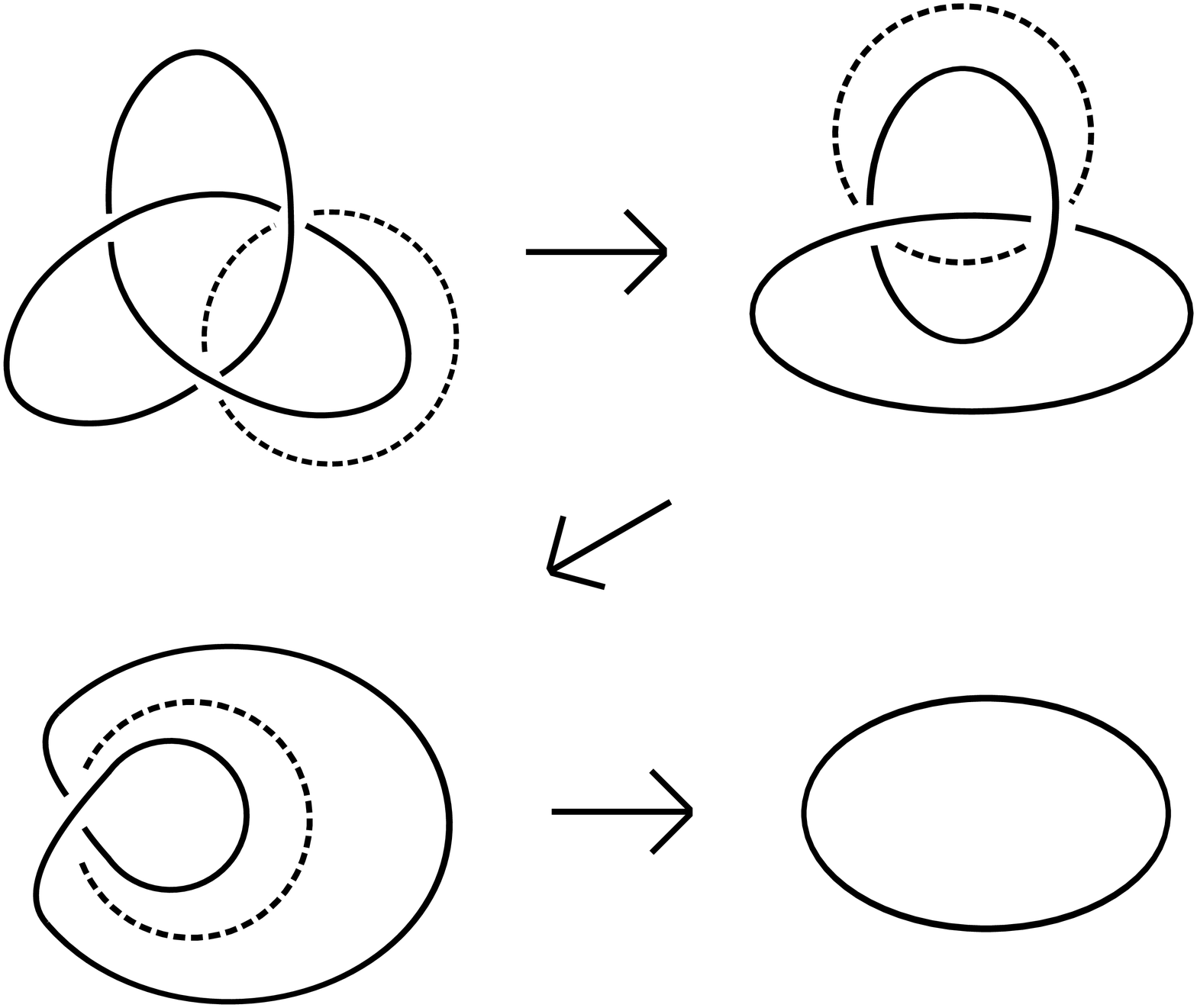} \\
\caption{Torefoil knot is in $D_{red}$}\label{trefoil}
\end{figure}

For example, the trefoil knot diagram is in $D_{red}$ (see Figure \ref{trefoil}). But the alternating diagram of the Borromean rings are not in $D_{red}$.

Let $\mathcal{L}_{red} = \{ L ; L \text{ is B-reducible} \}$, where $B$-reducible means there is some alternating link diagram $D_{L}$ of $L$ in $D_{red}$.

\subsection{Equivalence of $\mathcal{L}_{red}$ and $\mathcal{L}_{\overline{\rm{Brm}}}$ }

\begin{thm} \label{thmRedBrm}
$\mathcal{L}_{red} = \mathcal{L}_{\overline{\rm{Brm}}} = \{\text{ an alternating link } L; \rm{Brm} \nleq L\}$, where $\rm{Brm}$ is the Borromean rings.
\end{thm}

\begin{prf}
First, note some easy observations.
Let $L_{1}$ and $L_{2}$ be alternating links in $S^3$. If an alternating diagram $D_{L_{2}}$ of $L_{2}$ is given by reducing some $D_{L_{1}}$ of $L_{1}$ by (I)-move, then $D_{L_{1}} \subseteq D_{L_{2}}$ and $L_{1}=L_{2}$. On the otherhand, if $D_{L_{2}}$ is given by reducing $D_{L_{1}}$ by (II)-move, then $D_{L_{1}} \subseteq D_{L_{2}}$.

$\underline{\mathcal{L}_{red} \subset \mathcal{L}_{\overline{\rm{Brm}}}}$.
Assume that $L$ is a $B$-reducible alternating link which satisfies $\rm{Brm} \leq L$. We should conclude a contradiction.
By definition of $\mathcal{L}_{\overline{\rm{Brm}}}$, $\rm{Brm} \subseteq D_{L}$ for any minimal-crossing alternating link diagram $D_{L}$. Since $D_{L}$ is $B$-reducible, there is a sequence of finite disks ${B_{1},\cdots ,B_{n}}$ and there is some $m>0$ such that $\rm{Brm} \subseteq D_{L}(B_{1},\cdots ,B_{m})$ and $\rm{Brm} \nsubseteq D_{L}(B_{1},\cdots ,B_{m+1})$. So by the above observations, we can assume that $D_{L}$ satisfies $\rm{Brm} \nsubseteq D_{L}(B)$ for a reducible disk $B$ without loss of generality.

\begin{itemize}
\item When $B$ is $1$-reducible, $D_{L}$ is represented as a connected sum of two link diagrams (see Figure \ref{B1redcase}). Since the Borromean rings are irreducible, it is contained in one of the link diagrams. Then, $\rm{Brm} \subseteq D_{L}(B)$. This is a contradiction.
 
\item When $B$ is $2$-reducible, denote these two crossing points $c_{1}$ and $c_{2}$ and assume $c_{2}$ is smoothed by this operation (see Figure \ref{B2redcase}-(0)). By the assumption $\rm{Brm} \nsubseteq D_{L}(B)$, we should smooth some crossing points and they must contain $c_{1}$ or $c_{2}$. Otherwise, the Borromean rings contain this disk $B$ or $\rm{Brm} \subseteq D_{L}(B)$. These cases conclude contradictions. Thus, there remains five cases to smooth $c_{1}$ and $c_{2}$ (see Figure \ref{B2redcase}). But in each case, we can prove easily that $\rm{Brm} \subseteq D_{L}(B)$. Actually, we can prove similarly in the case of (2), (3), (4) and (5). In the case of (1), we observe that if there exists a disk $B$ in $\rm{Brm}$ whose boundary intersects one crossing point and two points of $D_{L}$, then the inside of $B$ is uniquely determined and we can prove $\rm{Brm} \subseteq D_{L}(B)$ (see Figure \ref{simplenbd}). This is contradiction.
\end{itemize}

\begin{figure}[h]

\includegraphics[width=4cm,clip]{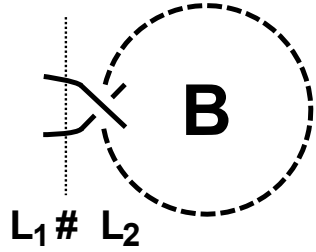} \\
\caption{1-reducible case}\label{B1redcase}
\end{figure}

\begin{figure}[h]

\includegraphics[width=8cm,clip]{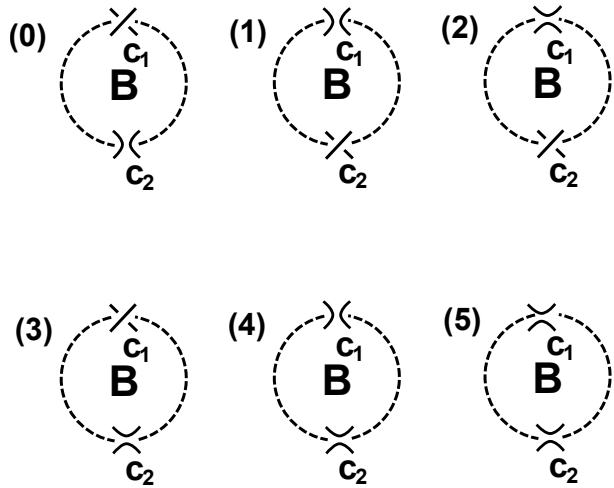} \\
\caption{2-reducible case}\label{B2redcase}
\end{figure}
\begin{figure}[h]

\includegraphics[width=4cm,clip]{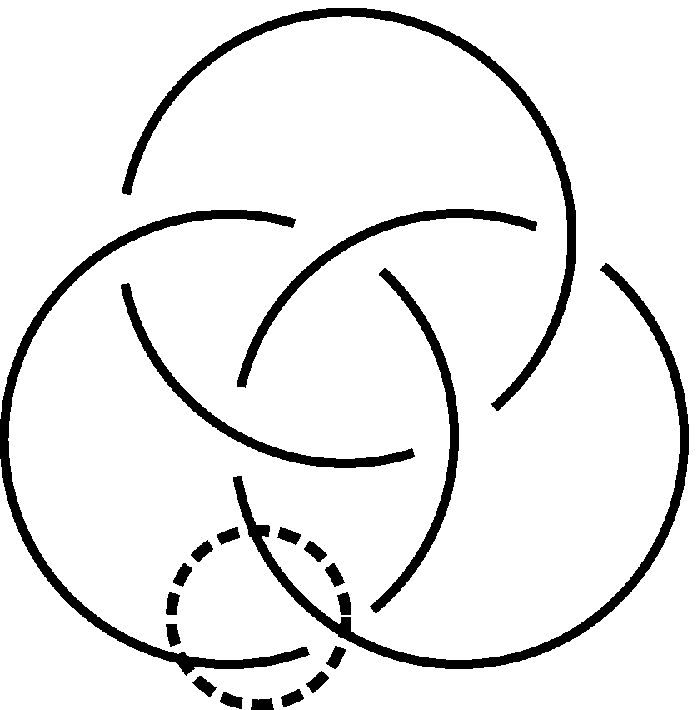} \\
\caption{}\label{simplenbd}
\end{figure}

$\underline{\mathcal{L}_{red} \supset \mathcal{L}_{\overline{\rm{Brm}}}}$.
Let $L$ be an alternating link which is not $B$-reducible. We  prove $\rm{Brm} \leq L$.
 
First, we can assume that an alternating link diagram $D_{L}$ of $L$ can not be represented as a disjoint union of alternating link diagrams. Otherwise, it is enough to consider the one of the components. We can also assume that under the above observation, $D_{L}$ satisfies the following condition (a).
\begin{itemize}
\item[(a):] $D_{L}$ does not admit any reducible disk and is not a disjoint union of the unknot diagrams (i.e., there exist some crossing points).
\end{itemize}

Then, it is enough to prove $\rm{Brm} \subseteq D_{L}$. Actually, if $\rm{Brm} \subseteq D_{L}(B_{1},\cdots, B_{n})$ for some reducible disks $(B_{1},\cdots, B_{n})$, then $\rm{Brm} \subseteq D_{L}$ and $\rm{Brm} \leq L$.

Next, let us call the closure of each component in $S^2 \setminus D_{L}$ a domain. Note that we can assume that each domain is wise. Otherwise, $D_{L}$ can be represented as a disjoint union of two alternating link diagrams.
Each domain $D$ has $k$ crossing points on its boundary(called $k$-gon). Note that $k \ge 3$ because $D_{L}$ does not admit any reducible disk. We  find $\rm{Brm}$ in $D_{L}$. 

Let $n_{k}$ be the number of $k$-gons $(k \ge 3)$ in $D_{L}$. Since the Euler number of $2$-sphere is two, we get the following equation by an easy computation.
\begin{equation}
n_{3} = 8 + \Sigma_{k \ge 5}(k-4)n_{k} \ge 8.
\end{equation} 
Thus, there are at least $8$ triangles in $D_{L}$. We  start with taking a triangle $D_{1}$. Let $\gamma_{1} = \partial D_{1}$. Since $D_{1}$ is a triangle then there are three polygons next to $D_{1}$. Let us denote these domains $D_{21}$, $D_{22}$ and $D_{23}$. We can prove that these domains satisfy the next three conditions.
\begin{itemize}
\item they are different domains,
\item they do not share their edges,  
\item they do not share their vertices other than the vertices of $D_{1}$. 
\end{itemize}
First, if $D_{21}$ and $D_{22}$ are the same domain, then there is a $1$-reducible disk $B$ (see Figure \ref{samedomain}).
Next, if $D_{21}$ and $D_{22}$ share their edges, then there is a disk $B$ whose boundary intersects with $D_{L}$ at three points (see Figure \ref{sharededge}). It is impossible. Lastly, if $D_{21}$ and $D_{22}$ share their vertices, then there is a $2$-reducible disk $B$. (see Figure \ref{sharedvertex}).

\begin{figure}[h]

\includegraphics[width=4cm,clip]{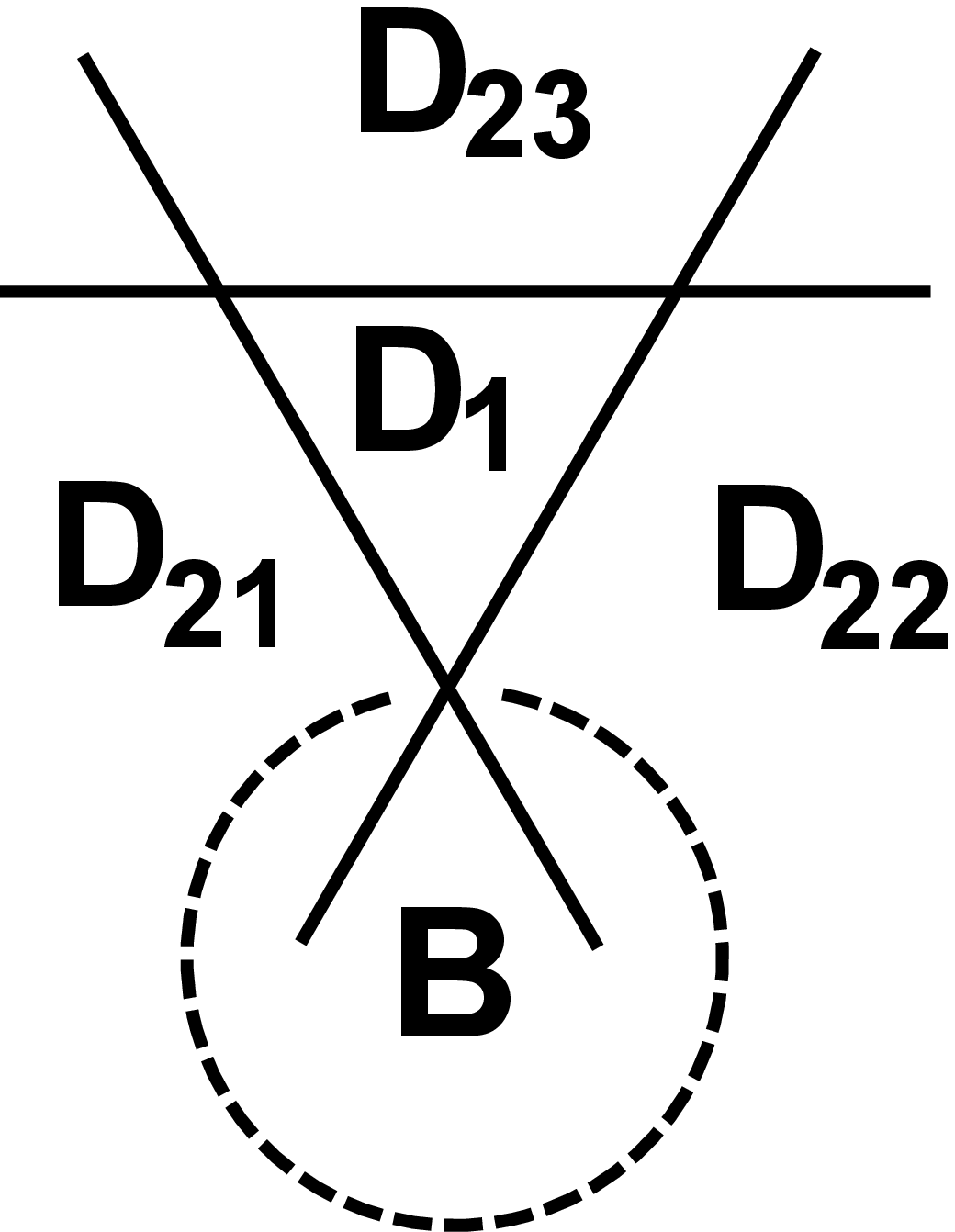} \\
\caption{}\label{samedomain}
\end{figure}

\begin{figure}[h]

\includegraphics[width=5cm,clip]{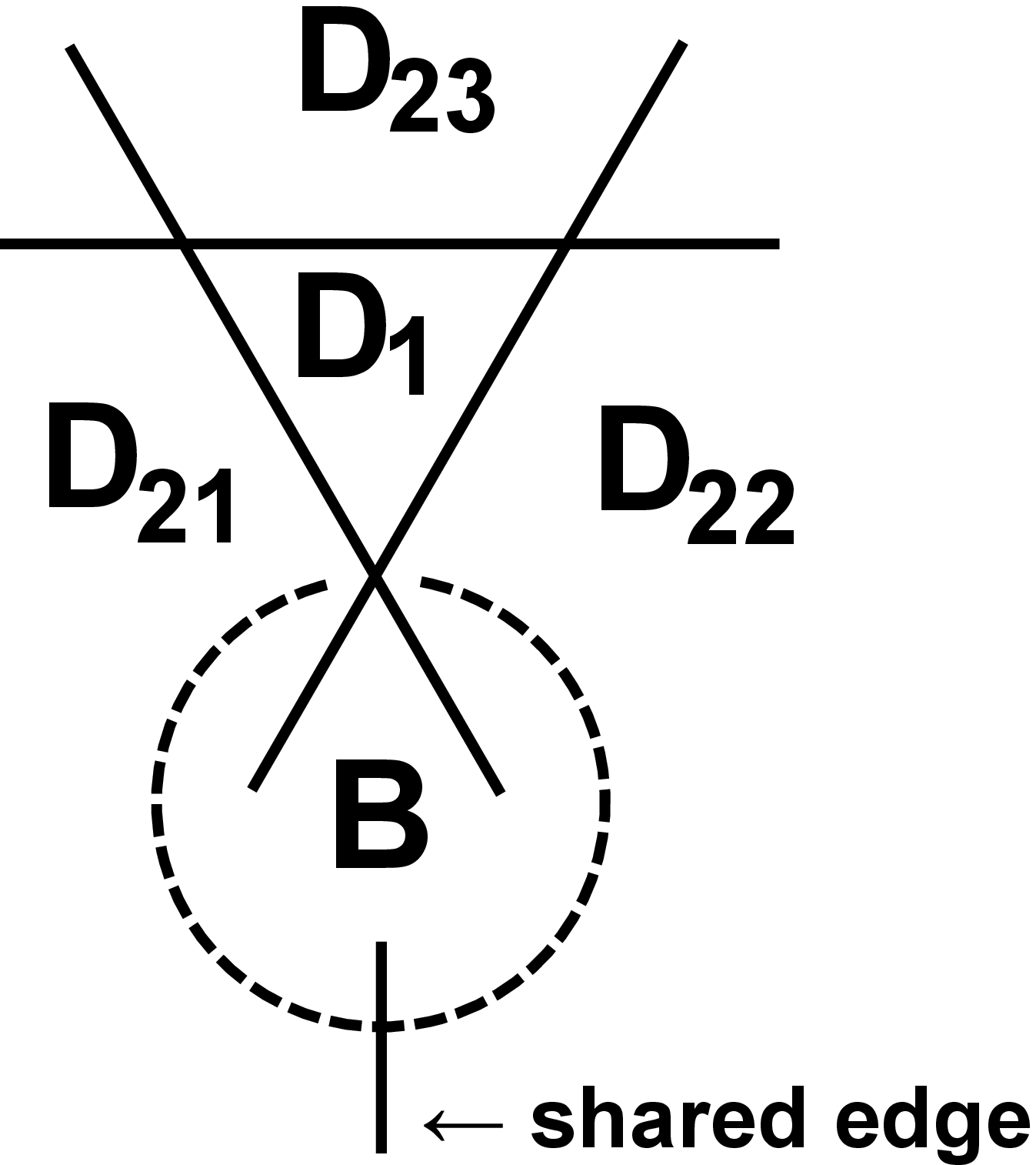} \\
\caption{}\label{sharededge}
\end{figure}

\begin{figure}[h]

\includegraphics[width=5cm,clip]{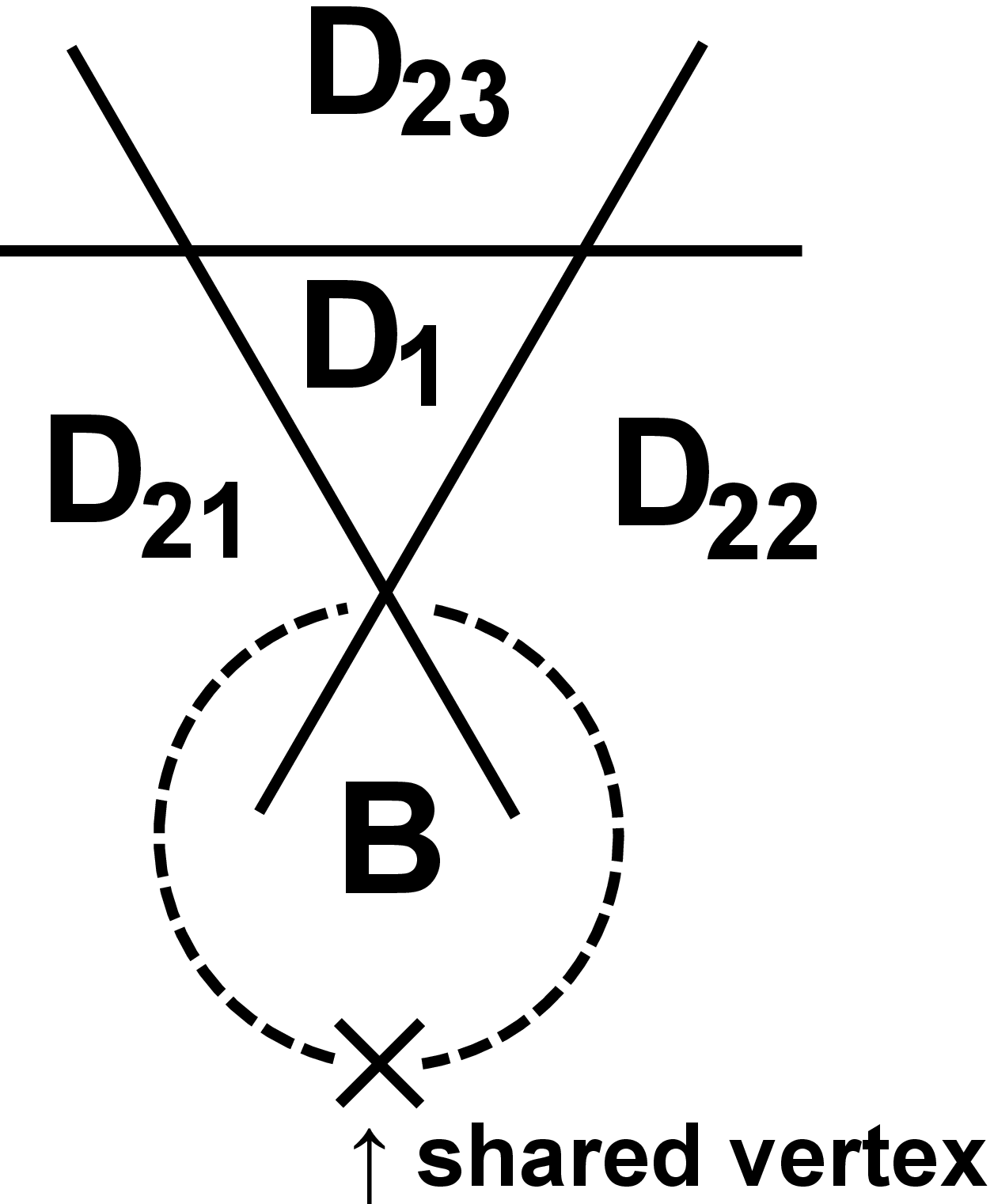} \\
\caption{}\label{sharedvertex}
\end{figure}

We put $D_{2} = D_{1} \cup (\cup_{i} D_{2i})$. We can regard $D_{2}$ as a polygon in $D_{L}$. Let $\gamma_{2} = \partial D_{2}$. Then, $\gamma_{2}$ is a simple closed curve. Next, there are $m$ domains next to $D_{2}$. Let us denote these domains $D_{31}, D_{32},\cdots, D_{3m}$. We can prove similarly that these domains do not share their edges. But it is possible that some of these domains coincide. We prepare the following definitions and two lemmas.
\begin{itemize}
\item[(b):] For an alternating link diagram $D_{L}$ satisfying (a) and a triangle $D_{1}$ as above, the domains  $D_{31}, D_{32},\cdots, D_{3m}$ are all disjoint.
\item[(c):]  For an alternating link diagram $D_{L}$ satisfying (a) and a triangle $D_{1}$ as above, some of these domains coincide. 
\item[(d):] For an alternating link diagram $D_{L}$, there are two different domains ${D^{*}}_{1}$ and ${D^{*}}_{2}$ which share their $l$-verticies ($l \ge 2$) such that $D_{L}$ admits just only $l(l-1)/2$ reducible disks. Moreover, there exist another one vertex or two vertices in $D_{L}$ when $l \ge 3$ or $l = 2$ respectively. (they are $2$-reducible disks (see Figure \ref{dtype}).　
\end{itemize}

\begin{figure}[h]

\includegraphics[width=5cm,clip]{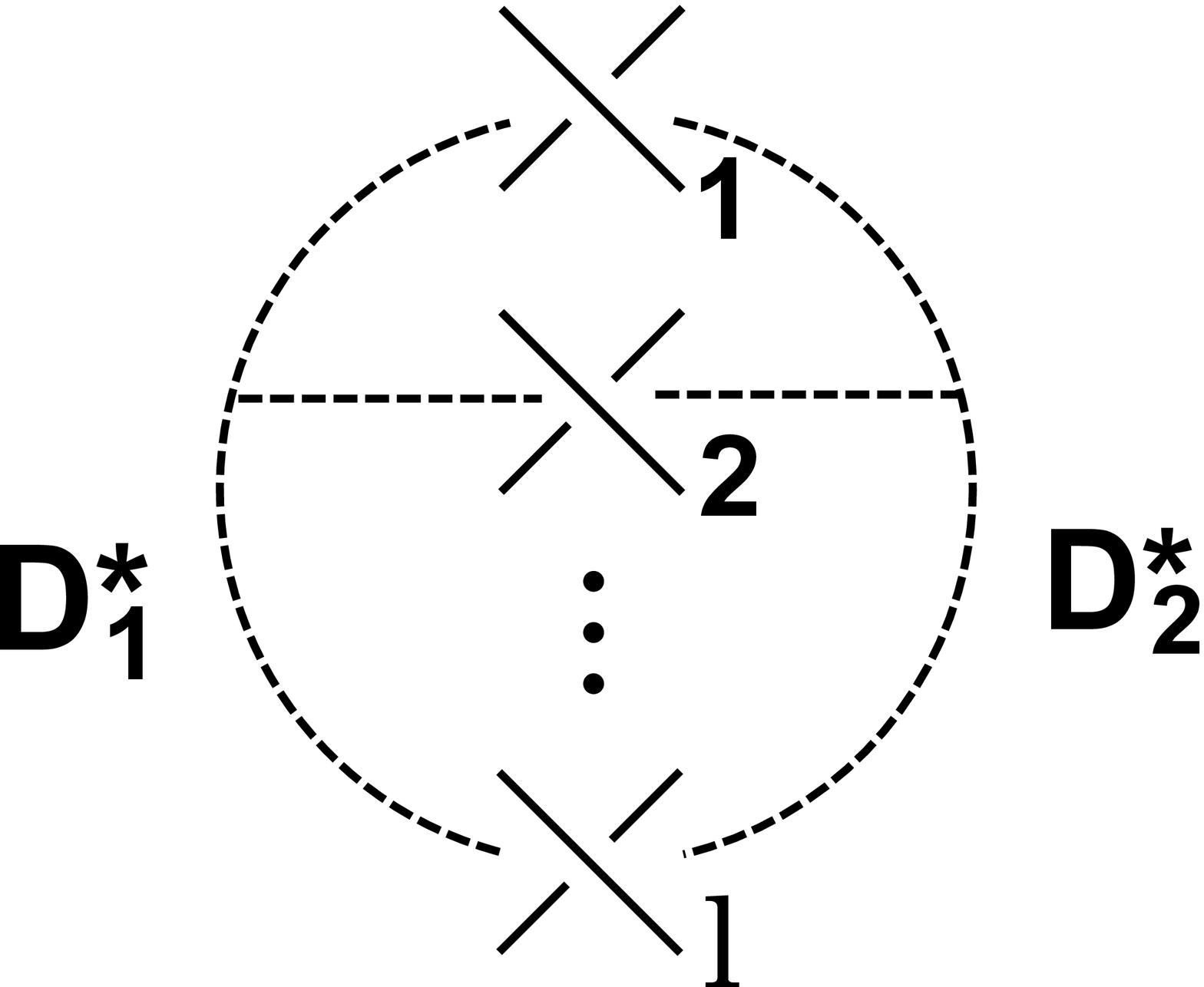} \\
\caption{}\label{dtype}
\end{figure}

\begin{lem} \label{lem1}
Let $D_{L}$ an alternating link diagram satisfies (a) and take a triangle $D_{1}$. Then, (b) or (c) occurs. Moreover, if (b) occurs, then we can find the Borromean rings in $D_{L}$, i.e., $\rm{Brm} \subseteq D_{L}$.
\end{lem}
\begin{lem} \label{lem2}
Under the condition (c) for $D_{L}$, we can find a new alternating link diagram $D'$ which satisfies (a) or (d) and $D' \subsetneq D_{L}$.
\end{lem}
\begin{lem} \label{lem3}
Under the condition (d) for $D_{L}$, we can find a new alternating link diagram $D'$ which satisfies (a) or (d) and $D' \subsetneq D_{L}$.
\end{lem}
Before proving these three lemmas, we first prove the proposition by using these lemmas.

Let $D_{L}$ be an alternating link diagram satisfying (a). Then, (b) or (c) occurs by Lemma \ref{lem1} If (b) occurs, we  finish the proof. If (c) occurs, then (a) or (d) occurs by Lemma \ref{lem2} Assume (d) occurs, then we can use Lemma \ref{lem3} just only finitely many times because the number of crossing points strictly decreases by these processes. So we will finally reach the condition (a). But we can also use Lemma \ref{lem2} just only finitely many times because the number of crossing points strictly decreases by this process. So we  reach the condition (b) by using Lemma \ref{lem1} Therefore, $\rm{Brm} \subseteq D_{L}$. 
\qed
\end{prf}

\begin{figure}[h]

\includegraphics[width=8cm,clip]{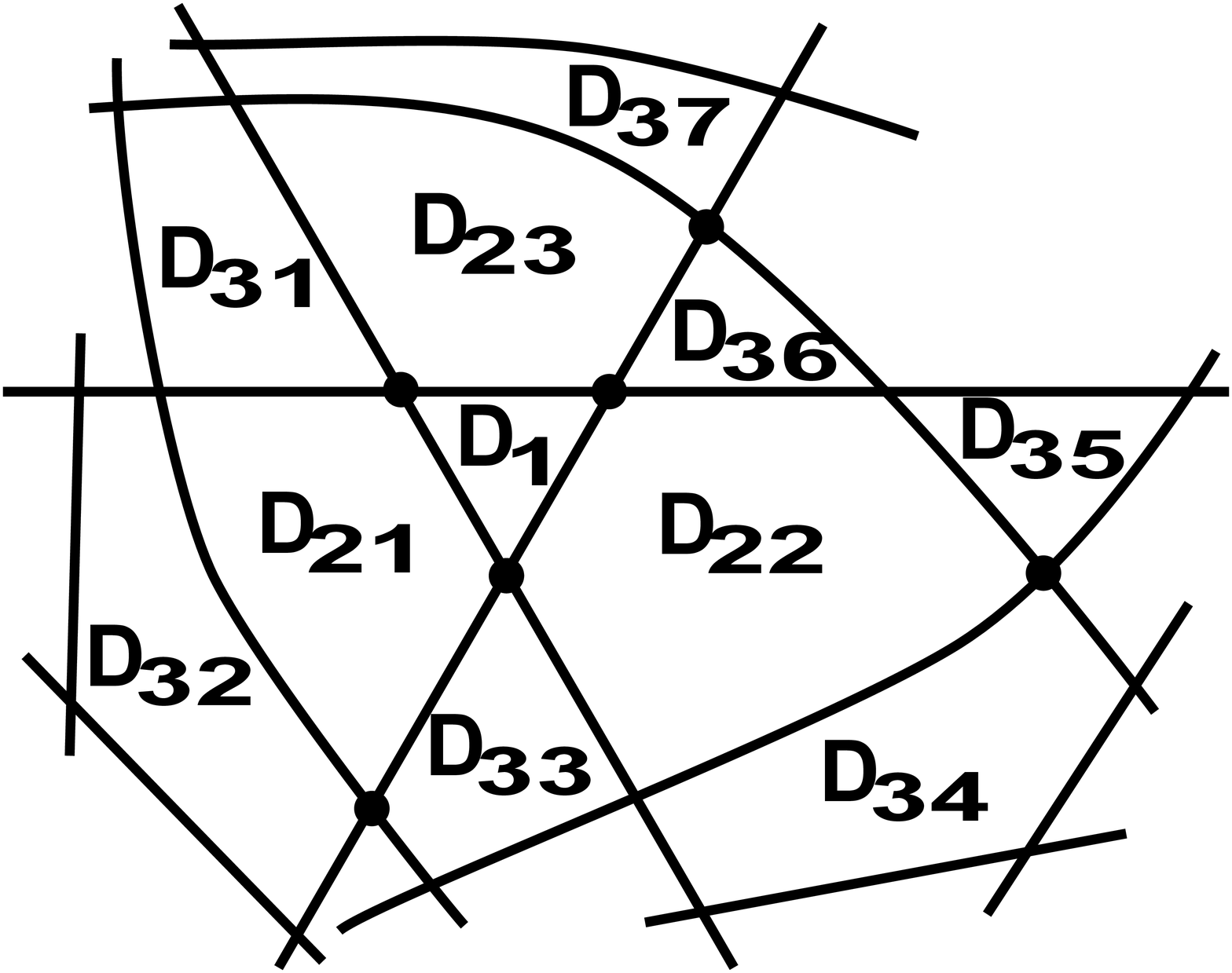} \\
\caption{}\label{findexample}
\end{figure}

\begin{figure}[h]

\includegraphics[width=8cm,clip]{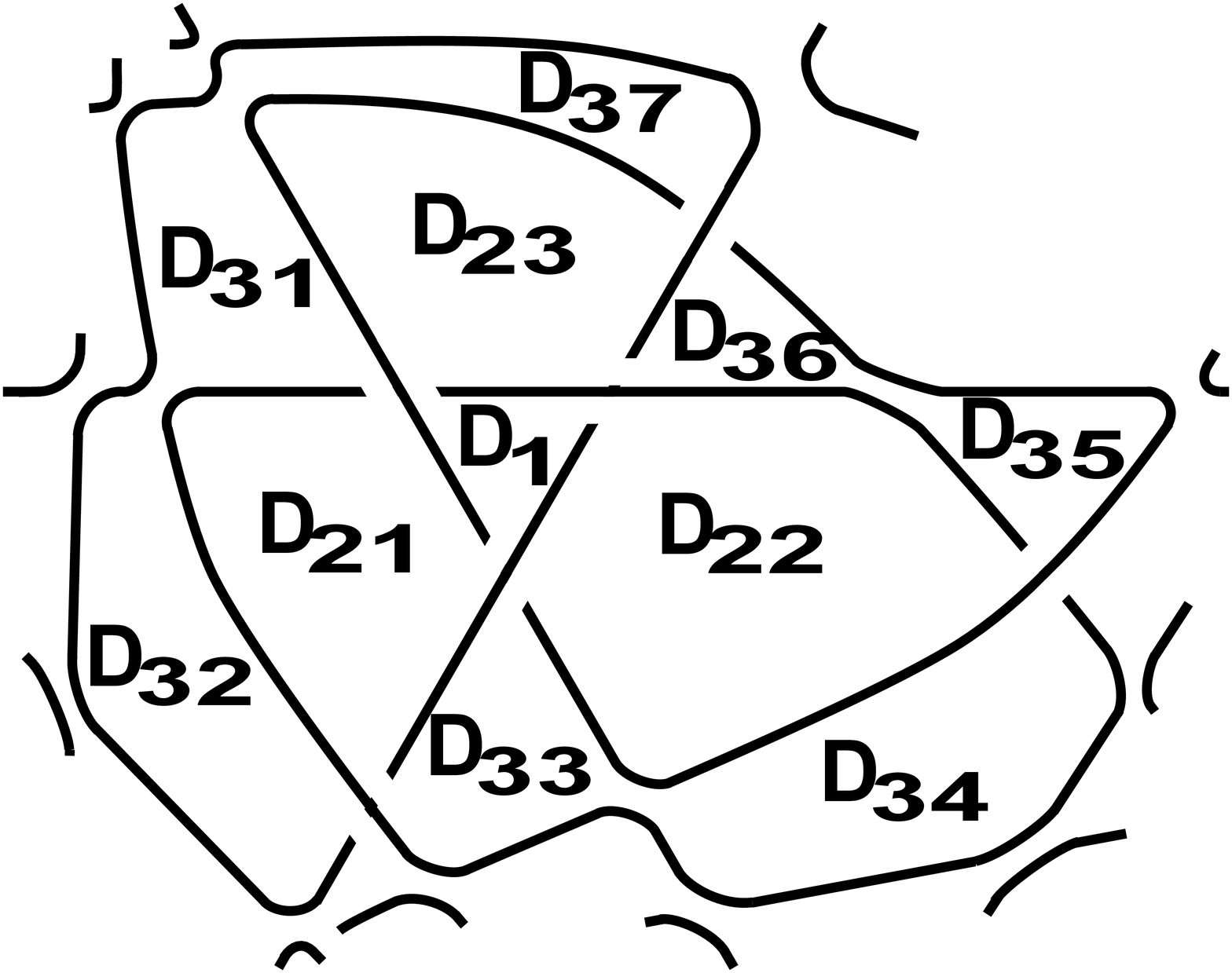} \\
\caption{}\label{findbrm}
\end{figure}

\begin{prf1}
The first part is trivial. If (b) holds, we put $D_{3} = D_{2} \cup (\cup_{i} D_{3i})$. Then, $D_{3}$ becomes a polygon in $D_{L}$ and let $\gamma_{3} = \partial D_{3}$. Thus, there are three curves $\gamma_{1}$, $\gamma_{2}$ and $\gamma_{3}$ and we can find $\rm{Brm}$ in $D_{L}$. Actually, we take the three vertices of the triangle $D_{1}$ and choose each vertex of $D_{21}$, $D_{22}$, $D_{23}$ respectively. With fixing these $6$ vertices, we smooth the other vertices along with three curves $\gamma_{1}$, $\gamma_{3}$ and $\gamma_{3}$ (see Figure \ref{findexample}). By this operation, we find $\rm{Brm}$ in $D_{L}$ (see Figure \ref{findbrm}).
\qed
\end{prf1}

\begin{figure}[h]

\includegraphics[width=8cm,clip]{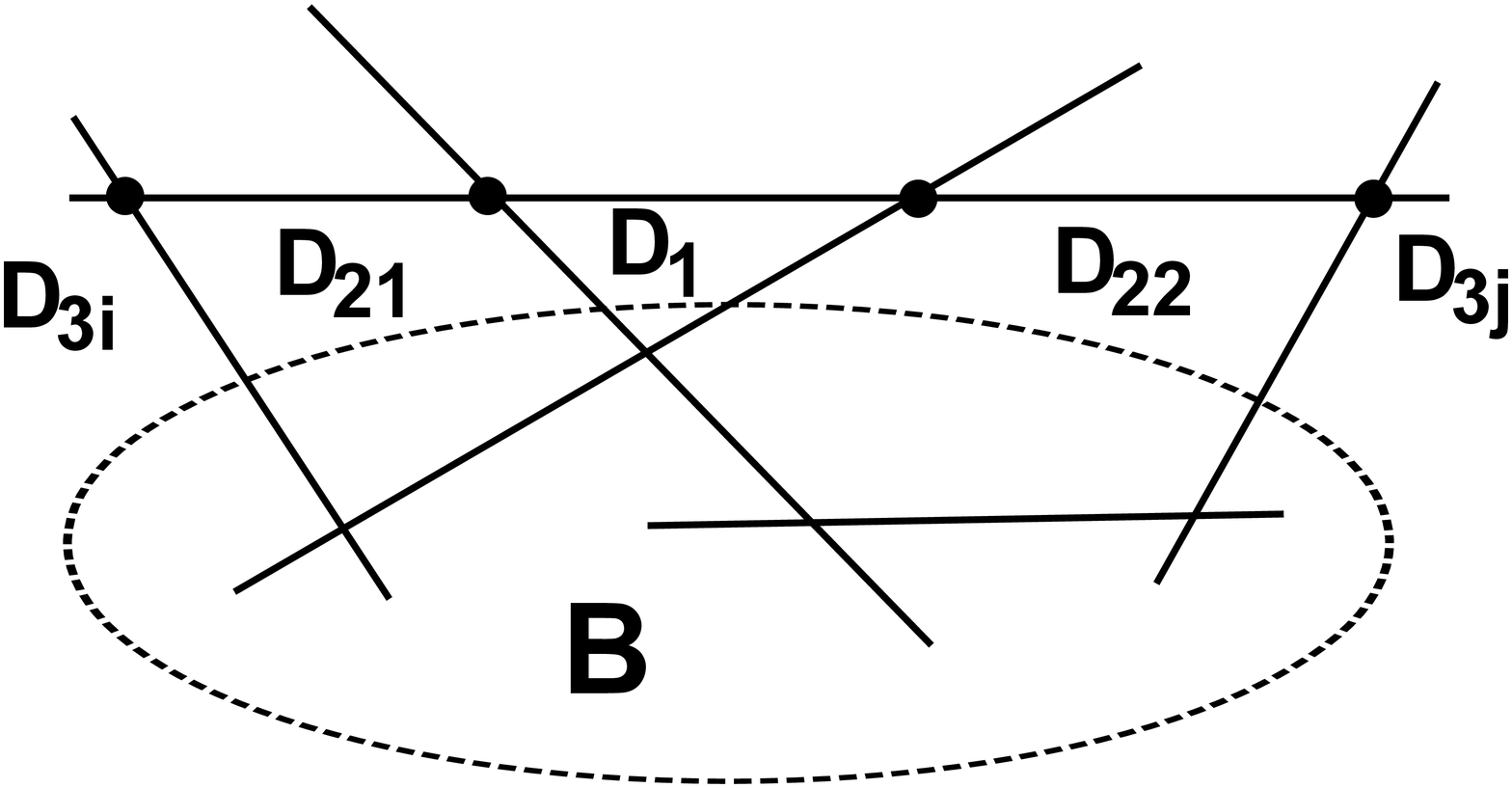} \\
\caption{}\label{inductivesamedomain}
\end{figure}

\begin{figure}[h]

\includegraphics[width=8cm,clip]{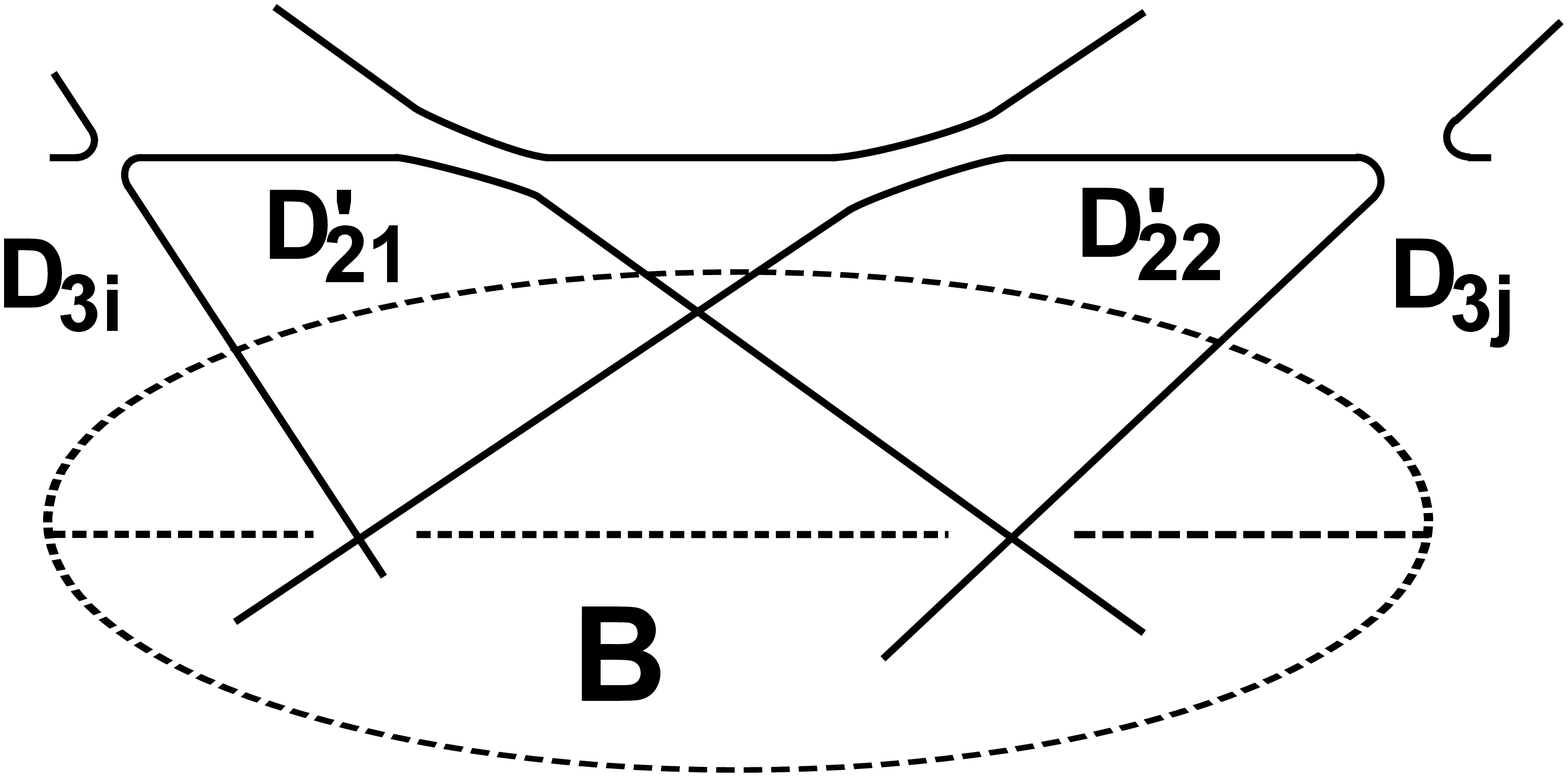} \\
\caption{}\label{doublebigon}
\end{figure}

\begin{figure}[h]

\includegraphics[width=8cm,clip]{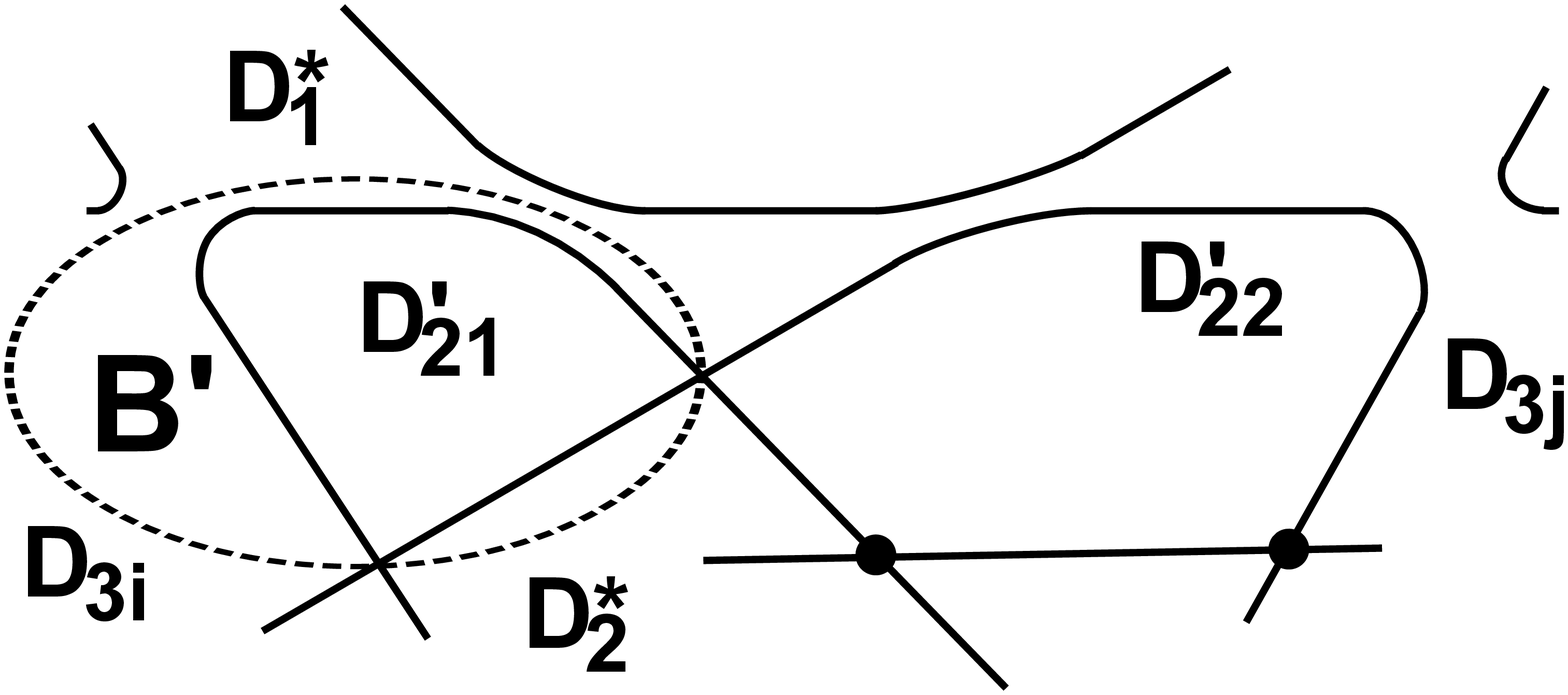} \\
\caption{}\label{onebigon}
\end{figure}

\begin{figure}[h]

\includegraphics[width=8cm,clip]{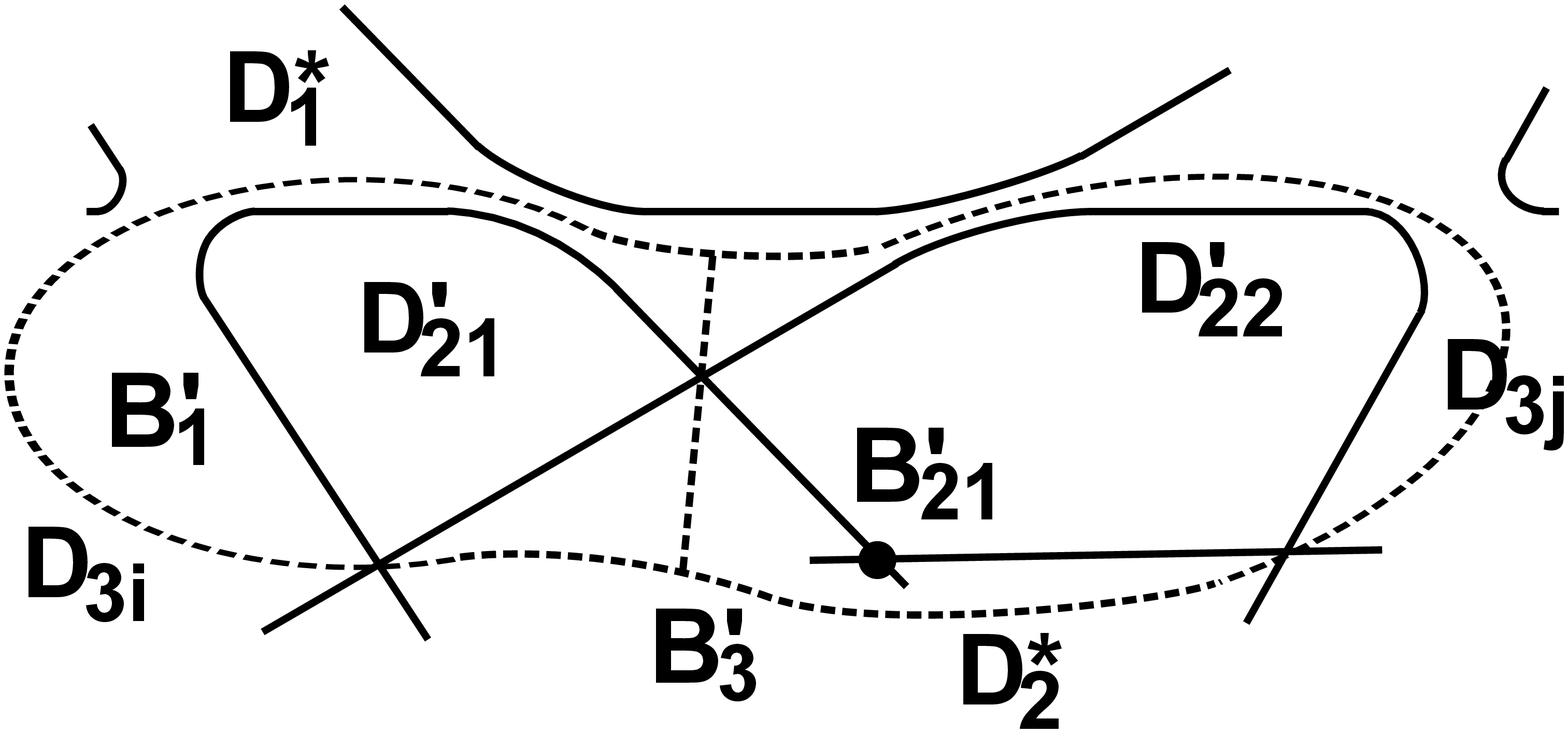} \\
\caption{}\label{threebigon}
\end{figure}

\begin{prf2}
We can assume that $D_{3i}$ and $D_{3j}$ coincide. When $D_{3i}$ and $D_{3j}$ are next to $D_{21}$, there is a $1$-reducible disk $B$. So this is a contradiction. Otherwise, we can assume that $D_{3i}$ is next to $D_{21}$ and $D_{3j}$ is next to $D_{22}$. There is a disk $B$ whose boundary intersects with $D_{L}$ at four points and intersects with $D_{3i}$, $D_{3j}$ and $D_{1}$ (see Figure \ref{doublebigon}). Take all vertices of $D_{21}$ and $D_{22}$ out of $B$. Then, smooth these vertices along the boundaries of $D_{21}$ and $D_{22}$. By this operation, we can find a new alternating link diagram $D'_{L} \subsetneq D_{L}$. Now we can assume one of new domains $D'_{21}$ and $D'_{22}$ has more than two vertices. Otherwise, there is a $2$-reducible disk in $D_{L}$ (see Figure \ref{inductivesamedomain}). Therefore, there remains three cases. 
\begin{itemize}
\item There is no reducible disk in $D'_{L}$. This is (a) for $D'_{L}$.
\item There is just one $2$-reducible disk $B'$ in $D'_{L}$ (see Figure \ref{onebigon}).
\item There are just three $2$-reducible disks $B'_{1}$, $B_{2}$ and $B'_{3}$ in $D'_{L}$  (see Figure \ref{threebigon}).
\end{itemize}
In the second and third cases, let ${D^{*}}_{1}$ and ${D^{*}}_{2}$ be the domains which intersect with $\partial B'$ or $\partial B'_{i}$. Therefore, $D'_{L}$ satisfies the condtion (a) or (d). 
\qed
\end{prf2}

\begin{figure}[h]

\includegraphics[width=4cm,clip]{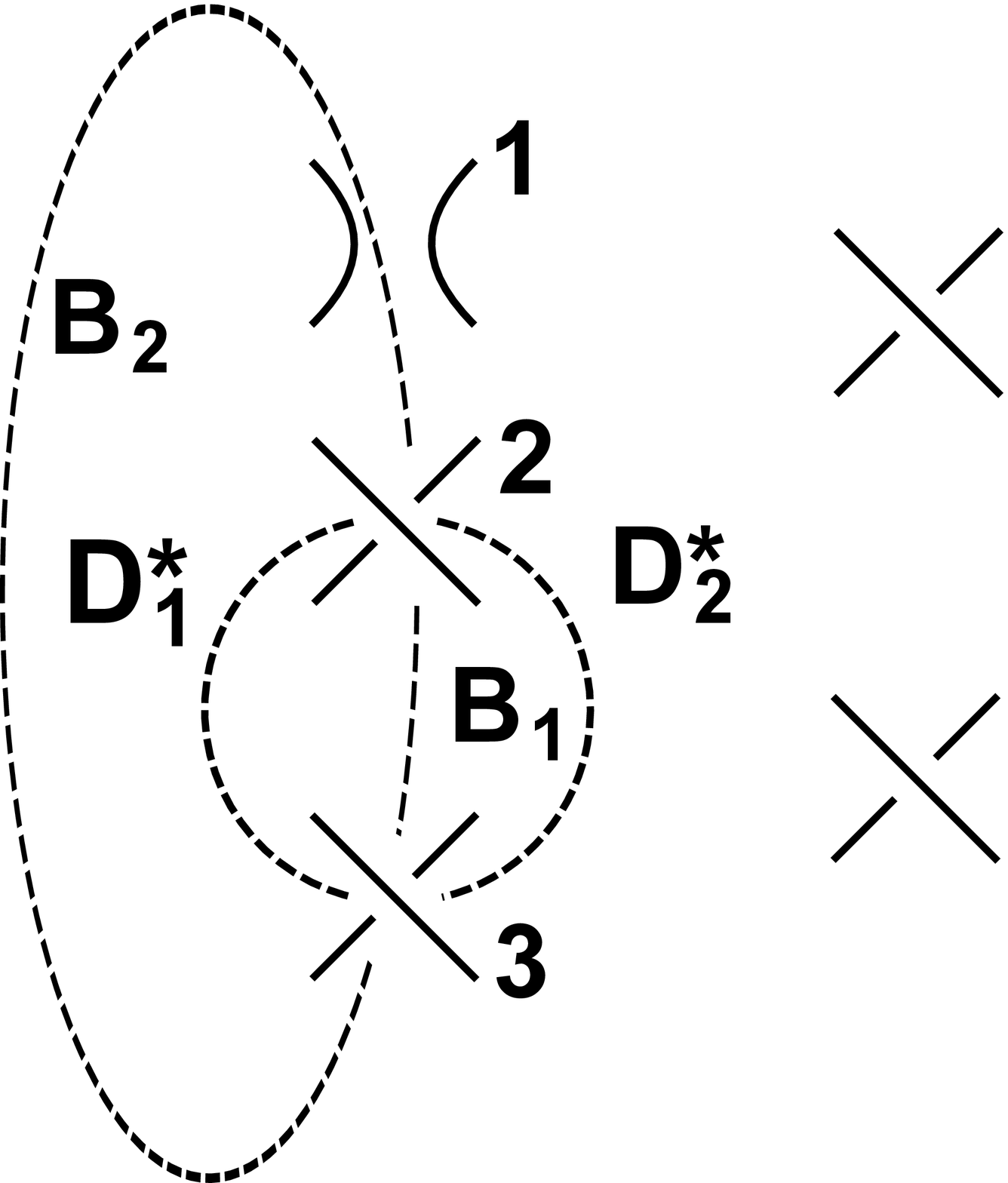} \\
\caption{}\label{d1}
\end{figure}

\begin{figure}[h]

\includegraphics[width=4cm,clip]{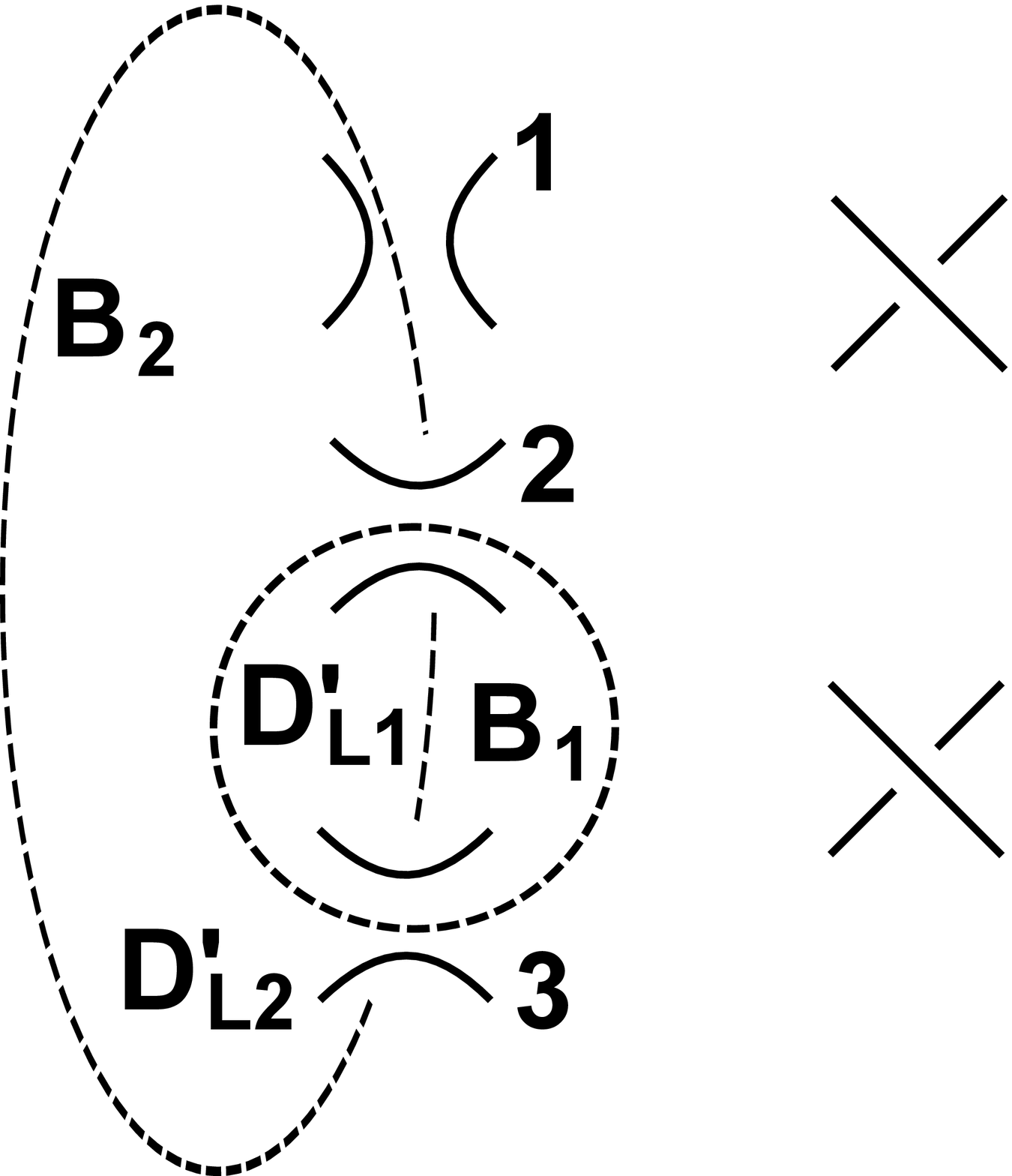} \\
\caption{}\label{d2}
\end{figure}

\begin{figure}[h]

\includegraphics[width=5cm,clip]{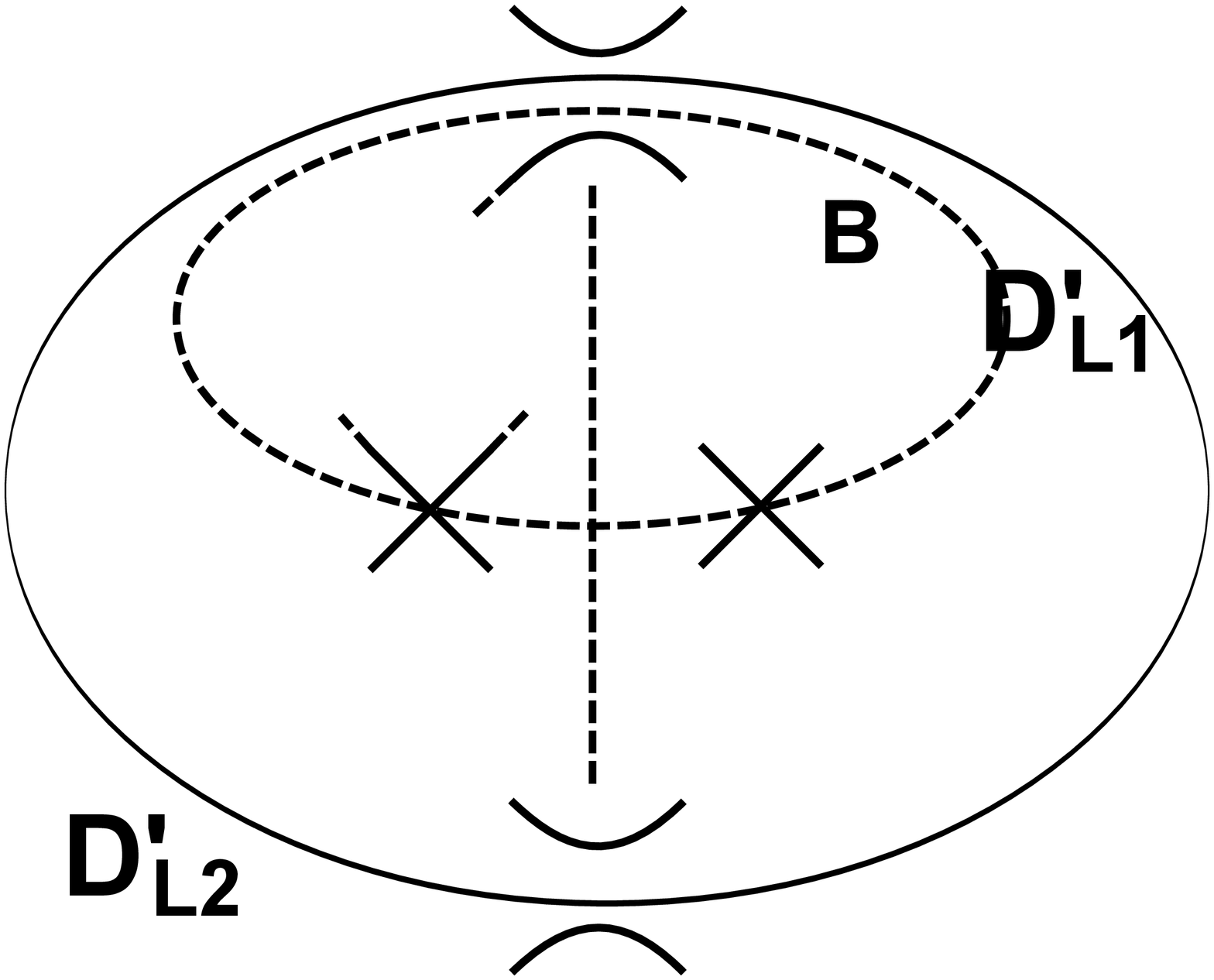} \\
\caption{}\label{d3}
\end{figure}

\begin{figure}[h]

\includegraphics[width=5cm,clip]{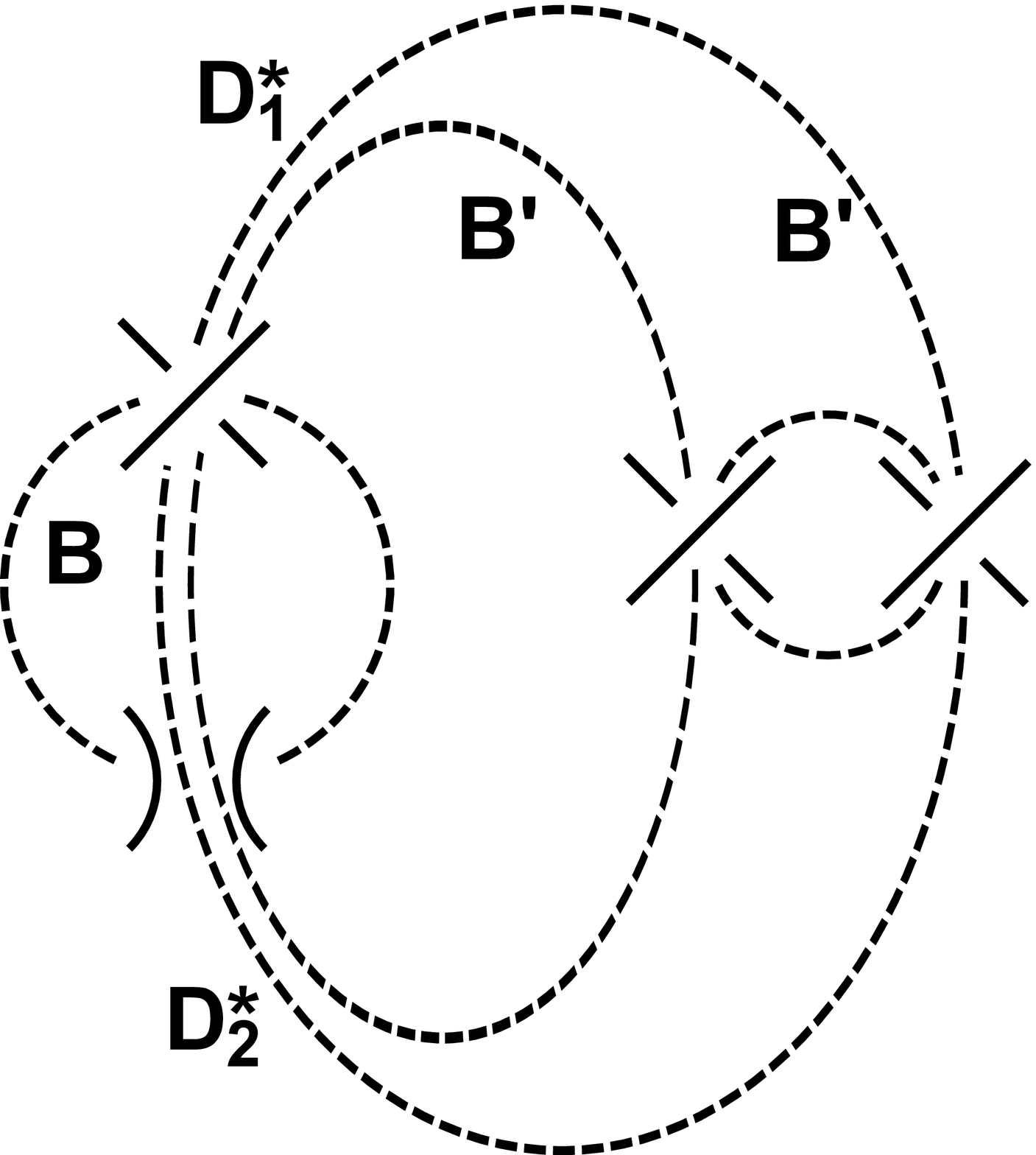} \\
\caption{}\label{d4}
\end{figure}

\begin{figure}[h]

\includegraphics[width=5cm,clip]{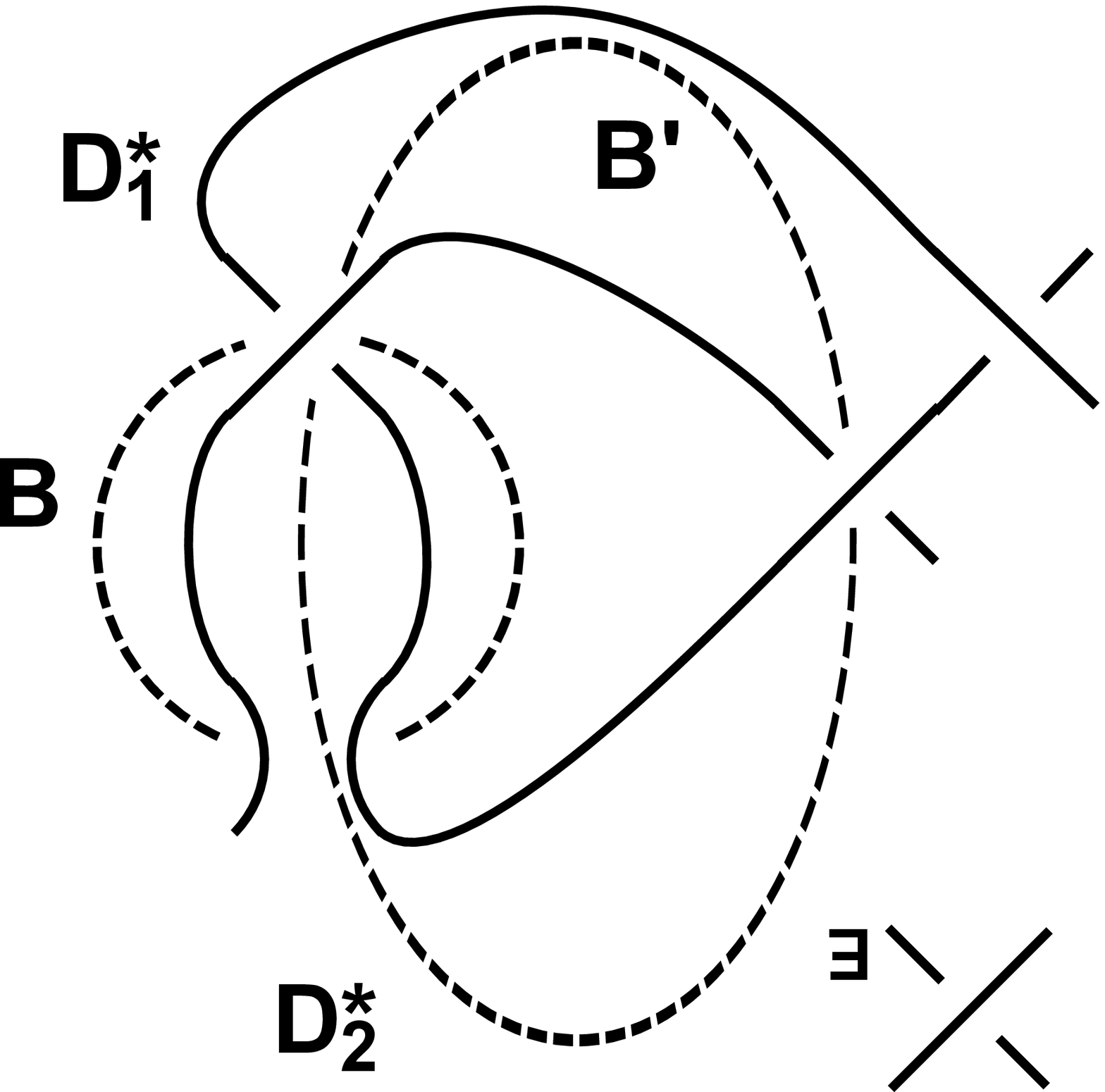} \\
\caption{}\label{morec}
\end{figure}

\begin{prf3}
First, we consider the case when $l > 3$. In this case, we can reduce $D_{L}$ by move-(II) until $l = 3$. This new diagram also satisfies (d). Next, we consider the case $l = 3$. In this case, we can reduce $D_{L}$ by move-(II) at $B_{1}$, but it happens that there is a new $2$-reducible disk $B_{2}$. Moreover, in this case, there is no other reducible disk because the boundary of such a reducible disk intersects with three crossing points in $D_{L}$ (see Figure \ref{d1}). We  smooth the rest two crossing points as in Figure \ref{d2}. Note that this new diagram can be represented as a disjoint union of alternating diagrams $D'_{L1}$ and $D'_{L2}$, i.e., inside and outside of $B_{1}$. We can prove at least one of this component $D'_{Li}$ satisfies the condition (a). Actually, if there exist $B$-reducible disks, $D_{L}$ can not become a diagram (see Figure \ref{d3}). Additionally there exist at least two crossing points in $D_{L}$, so at least one of $D'_{L1}$ and $D'_{L2}$ has a crossing point. Finally, we consider the case $l = 2$. In this case, we can reduce $D_{L}$ by move-(II). The new diagram $D'_{L}$ satisfies (a) or (d). Actually, if (a) is not satisfied, there exists only one reducible disk $B'$. If there exist other reducible disks, then $D_{L}$ has another reducible disks (see Figure \ref{d4}). Then, ${D^{*}}_{1}$ and ${D^{*}}_{2}$ are defined as the domains which intersect $\partial B'$ and we need at least another two vertices in $D'_{L}$. (see Figure \ref{morec}). Therefore, $D'_{L}$ satisfies (d).
\qed
\end{prf3}

%
\section{Alternatingly weighted trees and $B$-reducible alternating links}

\ In this section, we  first introduce another class of closed oriented three manifolds defined by surgeries along some links. After that, we claim that this class is also the same as the class of $B$-reducible alternating links. To prove this, we review the well-known correspondence between double branched coverings and Dehn surgeries (see \cite{Montesinos}). 

\subsection{Alternatingly weighted trees}

\begin{defn}
$(T,\sigma,w)$ is an \textit{alternatingly weighted tree} when the following three conditions hold.
\begin{itemize}
\item $T$ is a disjoint union of trees (i.e a disjoint union of simply connected, connected graphs). 
Let $V(T)$ denote the set of all vertices of $T$.
\item $\sigma:V(T) \rightarrow \{\pm 1\}$ is a map such that if two vertices $v_{1}$, $v_{2}$ are connected by an edge, then $\sigma(v_{1})=-\sigma(v_{2})$. 
\item $w:V(T) \rightarrow \{0, 1,\infty\}$ is a map.
\end{itemize}
\end{defn}  

\begin{figure}[h]

\includegraphics[width=9cm,clip]{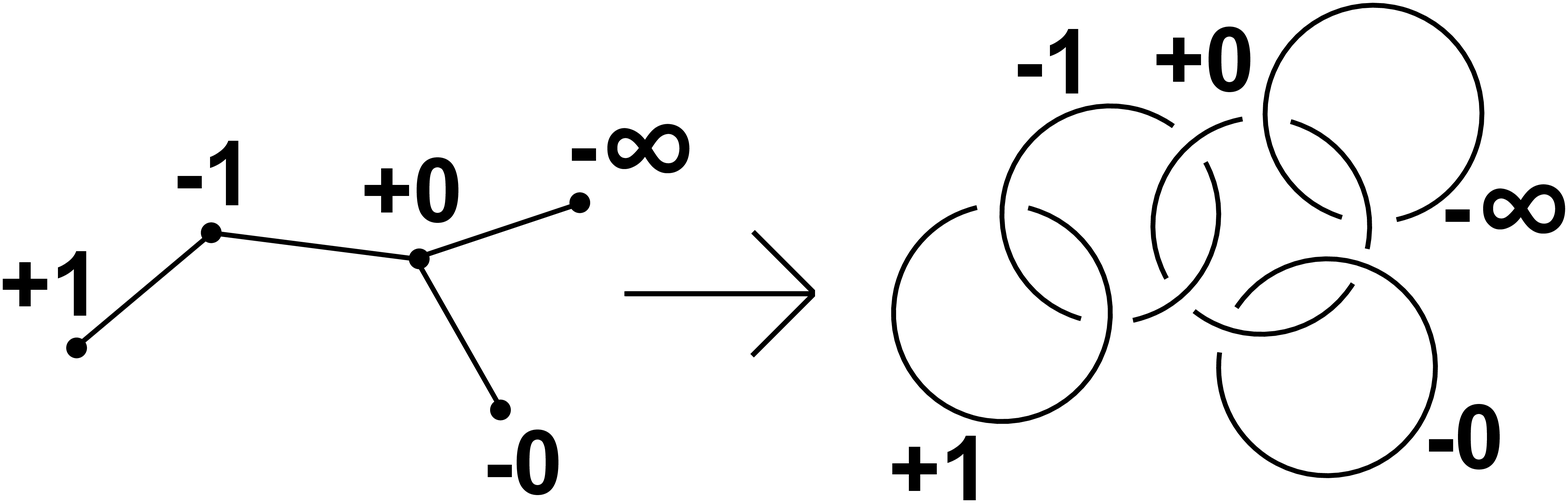} \\
\caption{}\label{tree1}
\end{figure}

Denote $\mathcal{T}$ the set of all alternatingly weighted trees.
For an alternatingly weighted tree $(T,\sigma,w)$, shortly $T$, we  define a three manifold $Y_{T}$ as follows.
First, we can take a realization of the tree $T$ in $\mathbb{R}^2 \subset S^3$. For each vertex $v$, we introduce the unknot in $S^3$. Next if two verteces in $T$ are connected by an edge, we link the corresponding two unknots with linking number $\pm 1$. Thus, we get a link $L_{T}$ in $S^3$. Then, we can get a new closed oriented three-manifold $Y_{T}$ by the surgery of $S^3$ along every unknot component of $L_{T}$ with the surgery coefficients $\sigma(v) w(v)$ (see Figure \ref{tree1})

This process gives a natural map $\mathcal{T} \rightarrow \mathcal{M}_{\mathcal{T}} = \{Y_{T} ; T \in \mathcal{T}\}/\text{homeo}$.

\begin{rem}
Note that we can also define the rational version of alternatingly weighred trees. That is, even if we replace the image $\{0, 1,\infty\}$ of $w$ by $\mathbb{Q}_{\ge 0} \cup \{\infty\}$, the induced manifolds are well-defined. In this case, we obtain rational surgeries of $S^{3}$ along links.
Moreover, the set of induced manifolds in the rational version are the same as $\mathcal{M}_{\mathcal{T}}$. Actually, we can represent a $\mathbb{Q}$-framed unknot by $\mathbb{Z}$-framed unknots, and we can also represent a $\mathbb{Z}$-framed unknot by $\{0, 1,\infty\}$-framed unknots by using continuous fraction expansions and slam-dunk operations, which is one of the Kirby calculus (see Figure \ref{qz1}).
\begin{figure}[h]

\includegraphics[width=9cm,clip]{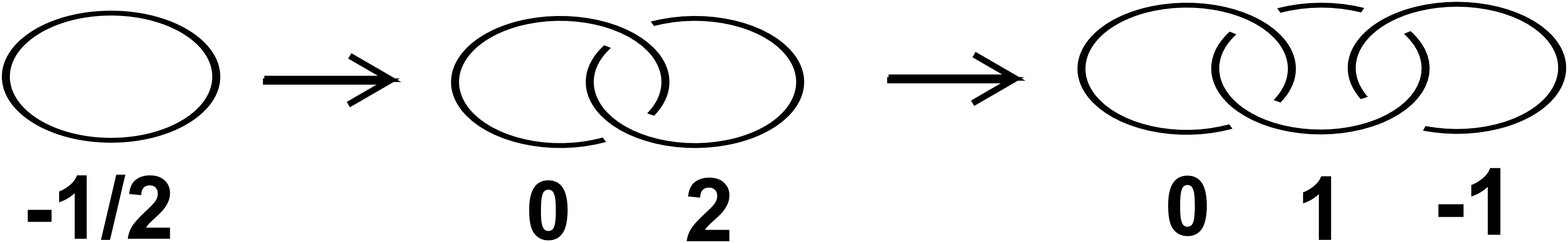} \\
\caption{}\label{qz1}
\end{figure}
\end{rem}



\begin{thm} \label{thmTred}
The set of the three manifolds $Y_{T}$ induced from alternatingly weighted trees $(T,\sigma,w)$ is equal to the set of the branched double coverings $\Sigma(L)$ of $S^{3}$ branced along $B$-reducible alternating link. That is,
$\mathcal{M}_{\mathcal{T}} = \mathcal{M}_{red} = \{\Sigma(L); L \text{ is in } \mathcal{L}_{red}\}$.
\end{thm}

\subsection{Montesinos Trick}

\begin{thm}\cite{Montesinos}
Let $M$ be a closed oriented three manifold that is obtained by doing surgery on a strongly-invertible link $L_{s}$ of $n$ components. Then, $M$ is a double branched covering of $S^3$ branched along a link $L_{d}$ of at most $n+1$ components. Conversely, every double branched covering of $S^3$ can be obtained in this fashion.
\end{thm}

A strongly-invertible link $L_{s}$ means a link with an orientation preserving invilution of $S^3$ which induces in each component of $L_{s}$ an involution with two fixed points. Without loss of generality, we can assume the involution is the axial symmetry $\phi$ with respect to $x$-axis.

We  sketch the method to get a new strongly-invertible link $L_{s}$ in $S^3$ from a link diagram $L_{d}$. It takes three steps.

\begin{enumerate}
\item Let $L_{d}$ be a connected link diagram (where connected diagram means a diagram which can not be written as a disjoint union of two diagrams.) For each crossing point $c$, we can take a small disk $B$ containing $c$ whose boundary intersects with $L_{d}$ at just four points. By smoothing each crossing, a new diagram has no crossing. Moreover, it is possible that the new diagram becomes the unknot diagram by smoothing suitably. (Of course this is not a unique way.) we assign the signature $+1$ or $-1$ to each disk by the following natural rules (see Figure \ref{step1-1} and \ref{step1-2}).

\begin{figure}[h]

\includegraphics[width=9cm,clip]{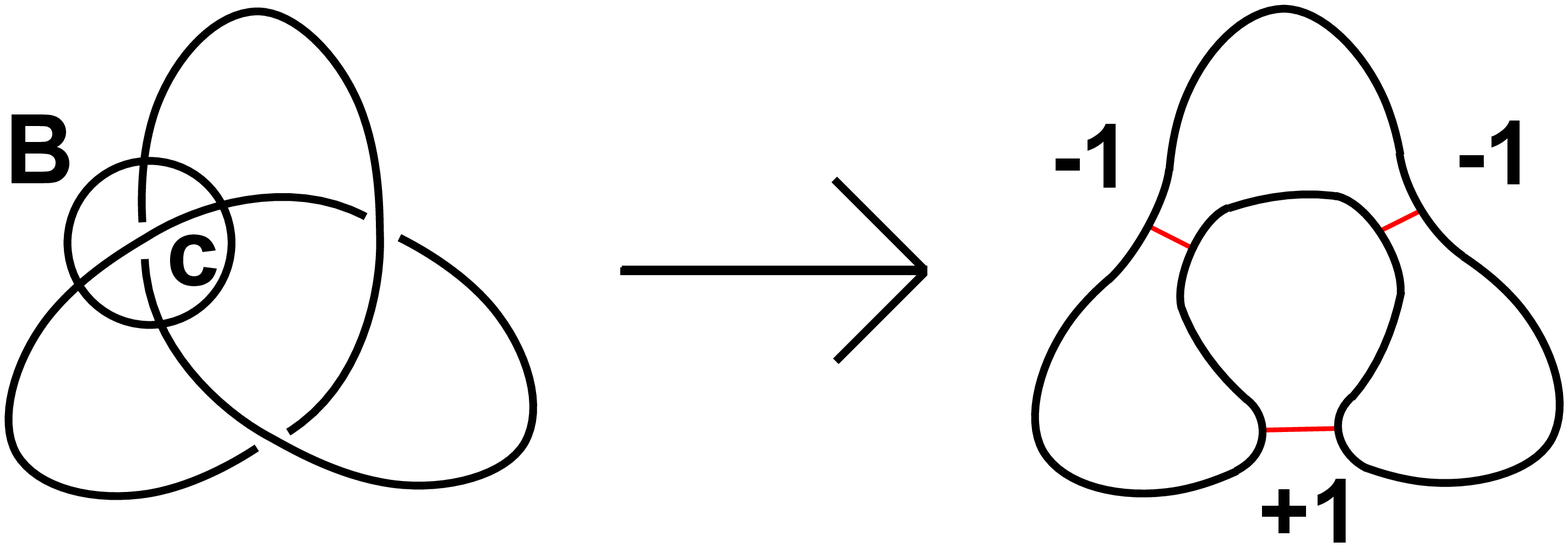} \\
\caption{}\label{step1-1}
\end{figure}
\begin{figure}[h]

\includegraphics[width=10cm,clip]{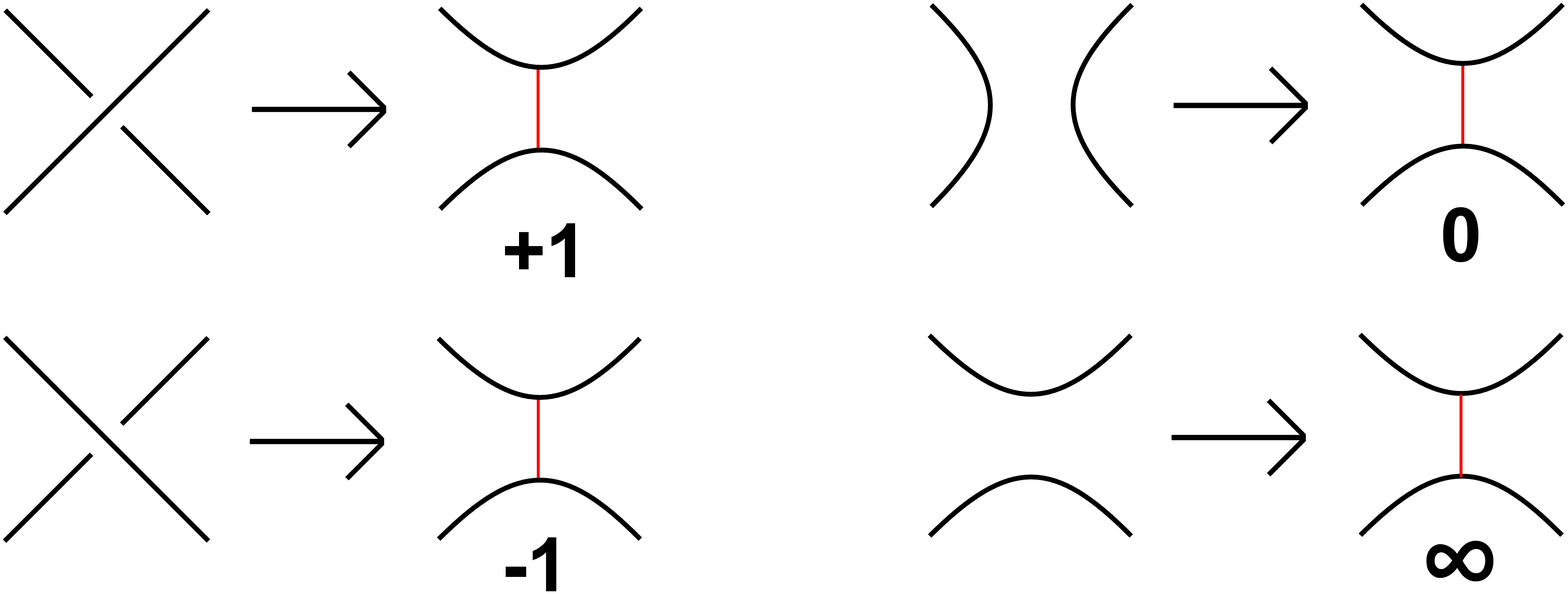} \\
\caption{}\label{step1-2}
\end{figure}

\item Since the new knot is just the unknot $U$, we can deform it by an isotopy to $x$-axis in $\mathbb{R}^3 \subset S^3=\mathbb{R}^3 \cup {\infty}$ by taking one marked point $p$ at $U \setminus (\text{disks})$ to the inifity. Let $\gamma_{B}$ be a trivial arc connecting the two arcs in each disk $B$. Then, $\{\gamma_{B}\}_{B}$ does not intersect each other (see Figure \ref{step2-1}).
\begin{figure}[h]

\includegraphics[width=10cm,clip]{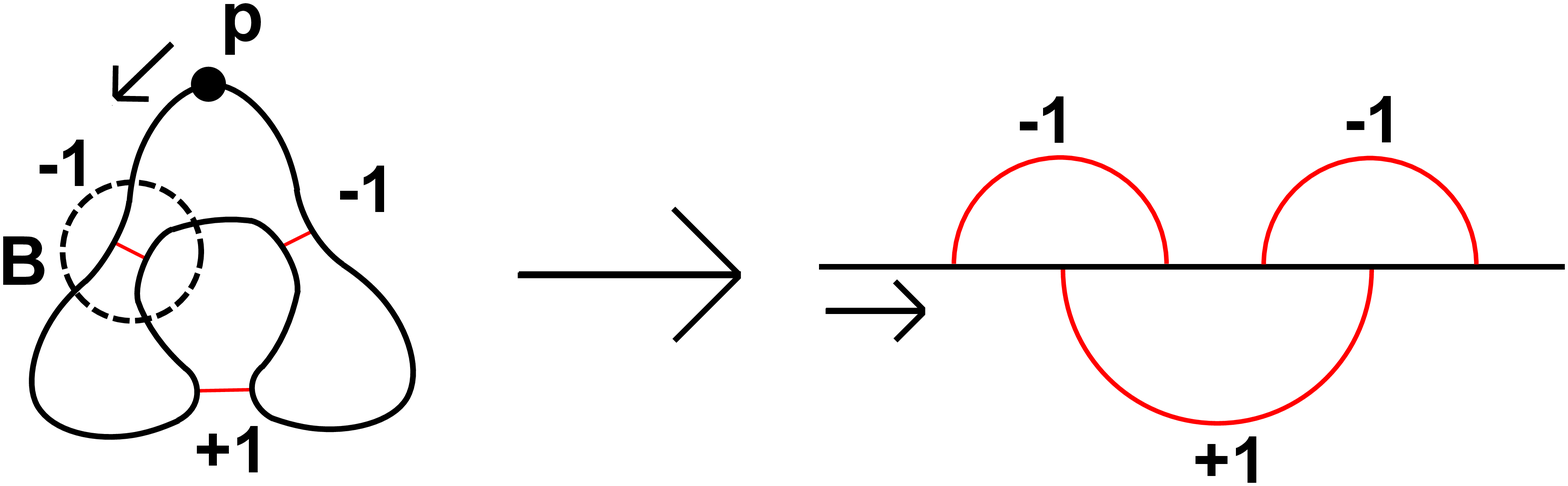} \\
\caption{}\label{step2-1}
\end{figure}

\item The double branched coverings of $S^3$ branched along the unknot is just $S^3$. Each arc $\gamma_{B}$ has its boundaries at $x$-axis. So a new link $L_{s}$ is defined as the double covering of $\gamma_{B}$ branched along the boundaries. By its definition, this new link $L_{s}$ is strongly-invertible. Moreover, each component of $L_{s}$ is the unknot (see \ref{step3-1}).

\begin{figure}[h]

\includegraphics[width=6cm,clip]{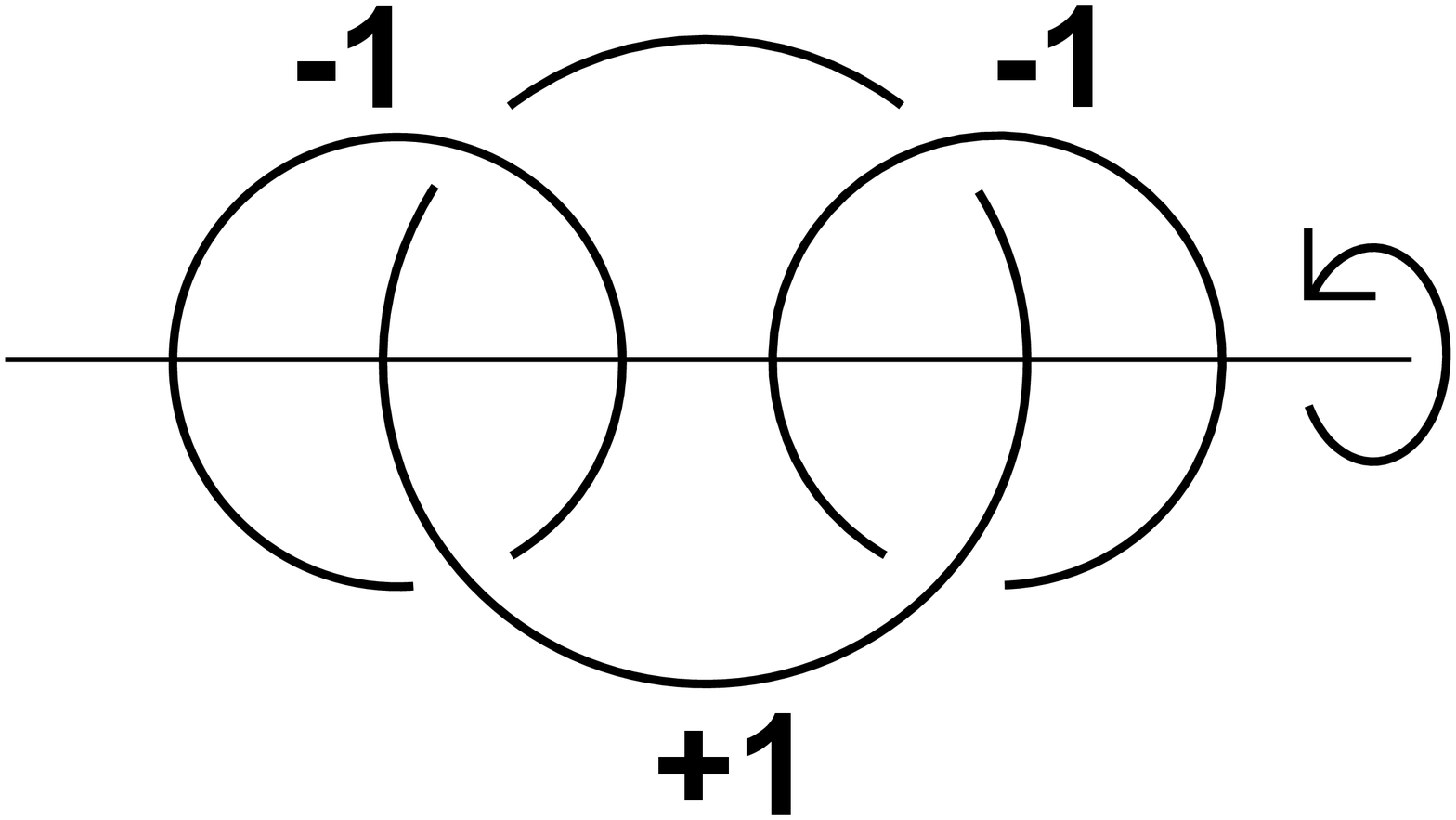} \\
\caption{}\label{step3-1}
\end{figure}

\end{enumerate}

\begin{defn}
For a strongly-invertible link $L_{s}$ and the involution $\phi$, the disjoint union $A_{L_{s}}$ of arcs in $\mathbb{R}^2$ whose boundary points are at $x$-axis is called a \textit{linear realization} of $L_{s}$ if $L_{s}$ corresponds to these arcs $A_{L_{s}}$ by the branched covering map $\phi$ from $S^3$ to $S^3$ branched along the $x$-axis $\cup \{\infty \}$. 
\end{defn}

Given a strongly-invertible link $L_{s}$ and its linear realization, we can get a link daigram $L_{d}$ by  reverse operations.
We call $L_{d}$ an $M$-induced link diagram of $L_{s}$.
Next, we  sketch the proof of the fact that the above method gives the equation $\Sigma(L_{d}) = S^3(L_{s})$.

\subsection{Proof of Theorem \ref{thmTred}}

We introduce the following four invertible operations to construct $T$ inductively.

\begin{figure}[h]

\includegraphics[width=6cm,clip]{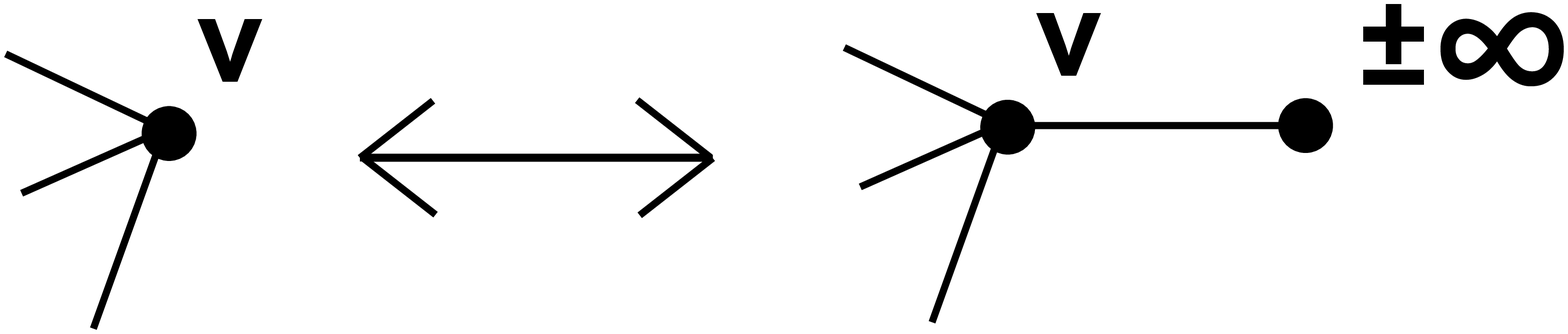} \\
\caption{}\label{1operation}
\end{figure}

\begin{figure}[h]

\includegraphics[width=6cm,clip]{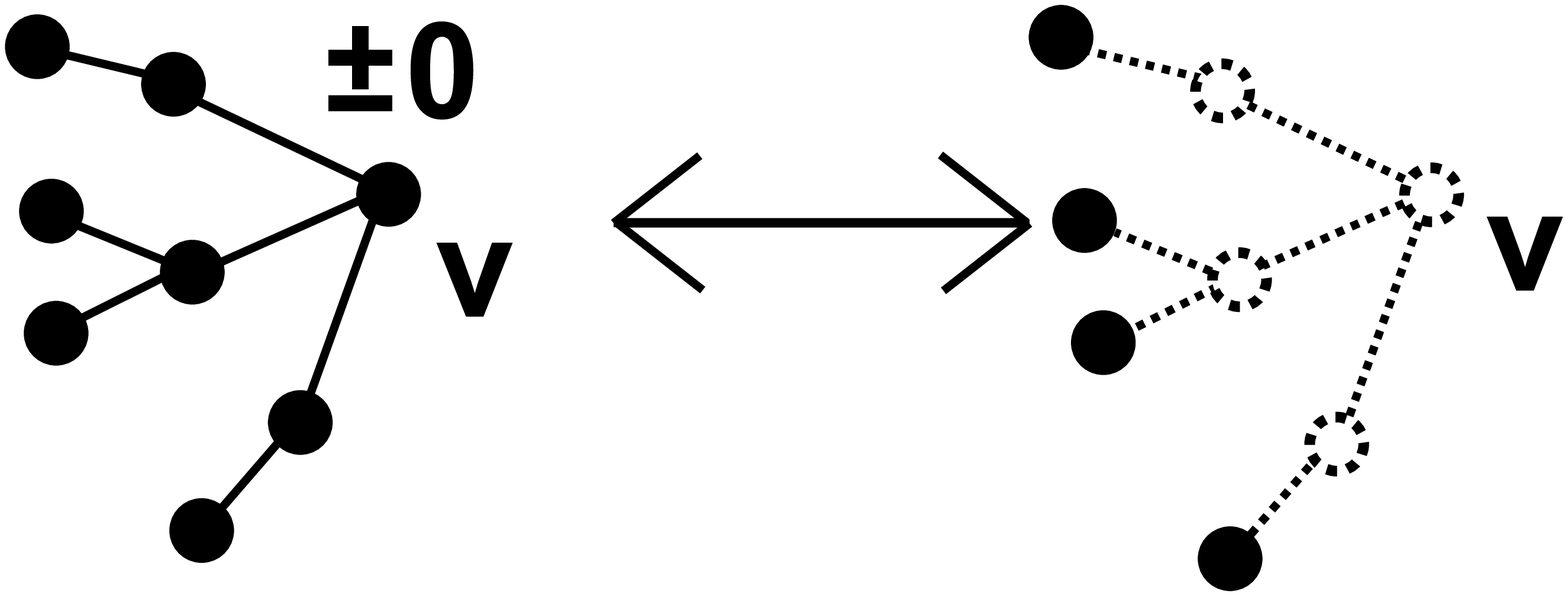} \\
\caption{}\label{2operation}
\end{figure}

\begin{figure}[h]

\includegraphics[width=6cm,clip]{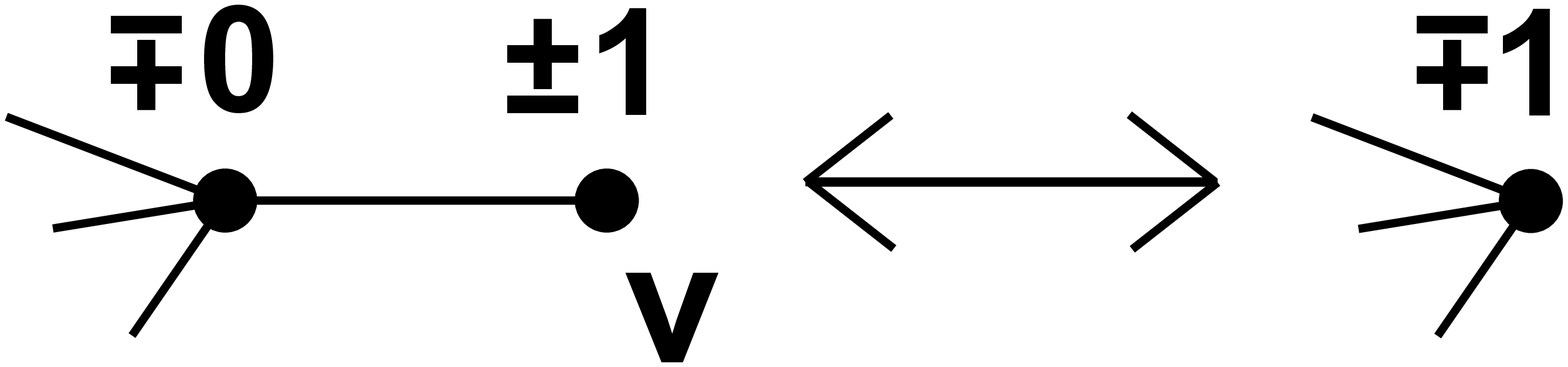} \\
\caption{}\label{3operation}
\end{figure}

\begin{figure}[h]

\includegraphics[width=6cm,clip]{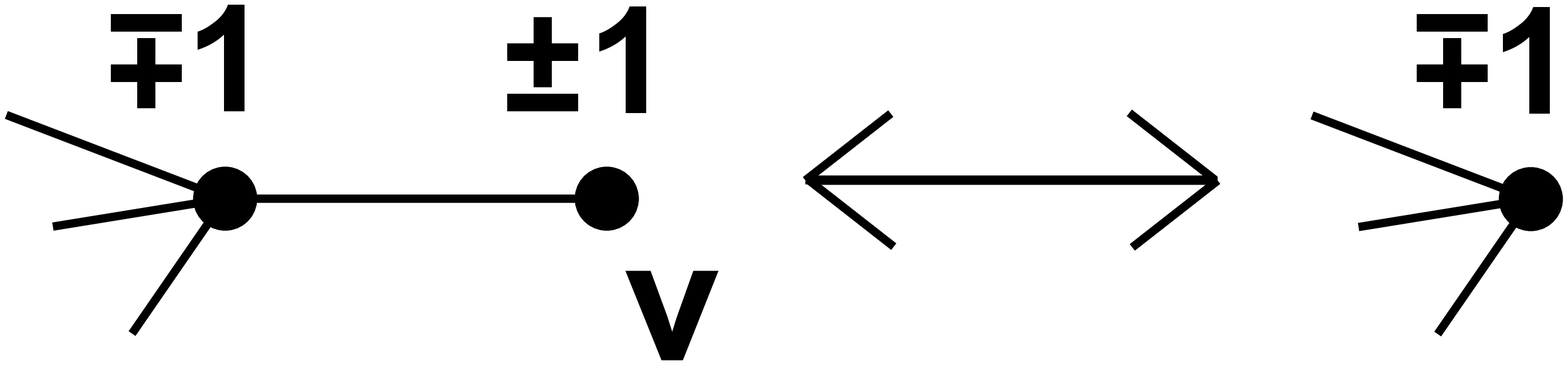} \\
\caption{}\label{4operation}
\end{figure}

\begin{enumerate}
\item For a vertex $v$, introduce a new vertex with weight $\mp \infty$ and connect it to $v$ (see Figure \ref{1operation}).
\item For a univalent vertex $v$ with weight $\pm 0$, remove the vertex, the next vertices and connected edges (see Figure \ref{2operation}).
\item For a univalent vertex $v$ with weight $\pm 1$, if the next vertex has its weight $\mp 0$, remove the univalent vertex and the connected edges, and change the weight of the next vertex into $\mp 1$ (see Figure \ref{3operation}).
\item For a univalent vertex $v$ with weight $\pm 1$, if the next vertex has its weight $\mp 1$, remove the univalent vertex and the connected edges (see Figure \ref{4operation}).
\end{enumerate}


\begin{rem}
The operations (1)-(3) do not change the induced three manifold $Y_{T}$. But operation (4) may change the induced three-manifold.
Moreover, each $T \in \mathcal{T}$ can be constructed from disjoint union of points with weight $\pm 0$ by using these operations in finitely many times.
Actually, for each $T$, use (1) to remove vertices with weight $\pm \infty$. Then, we can assume $T$ is connected. Using (2), (3) or (4), we can decrease the number of vertices of $T$. As a result, $T$ may be assumed to have only one vertex with weight $0$ or $\pm 1$. However, a point with weight $\pm 1$ is vanished by using (3) and (2).
\end{rem}

We start to prove Theorem \ref{thmTred}.

\begin{prf}

$\underline{\mathcal{M}_{\mathcal{T}}\subset \mathcal{M}_{red}}$.
We prove the next claim by induction on the maximal number of vertices of each connected component of $T$. Denote it by $|T|$. 
\begin{clm}
Let $T \in \mathcal{T}$ and $Y_{T} \in \mathcal{M}_{\mathcal{T}}$. Then, for any linear realization $A_{T}$, the $M$-induced link diagram $L_{T}$ satisfies $L_{T} \in D_{red}$ and $\Sigma(L_{T}) = Y_{T}$ 
\end{clm}
If $|T|=1$, $T$ becomes a disjoint union of points with weight $0$ or $\pm 1$. In this case, the linear realization of $T$ is unique and whose $M$-induced link $L_{T}$ becomes the unknot (see Figure \ref{t1case}).

\begin{figure}[h]

\includegraphics[width=10cm,clip]{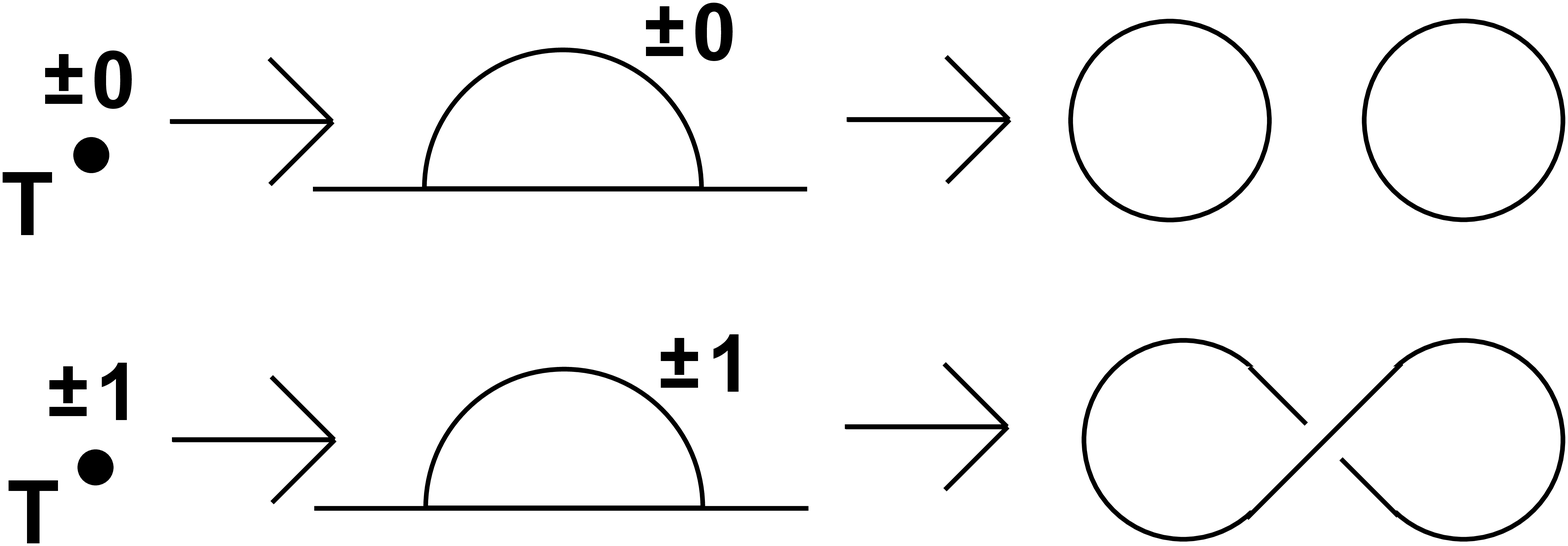} \\
\caption{}\label{t1case}
\end{figure}

Next, assume that the proposition holds when $|T| \leq n$. Take $T \in \mathcal{T}$ with $|T| = n+1$.
Then, Remark 4.1 tells us that $T$ can be changed into $T'$ with less vertices than $T$ by a operation (1), (2), (3) or (4). We consider case-by-case.

\begin{enumerate}
\item For any linear realization $A_{T}$, the natural linear realization $A_{T'}$ is induced by ignoring $\infty$ arc. Let $L_{T} = L_{T'}$. Then, $L_{T} \in D_{red}$ and $\Sigma(L_{T}) = Y_{T'} = Y_{T}$ by induction.
\item In this case, an arbitrary linear realization of $T$ looks  as in Figure \ref{2-1-1} and \ref{2-1-3}. We express some collections of arcs by numbers (i),..,(iv). Note that there is no arc connecting (i) and (v) with (ii), and (iii) with (iv). So we express this situation by $\times$. Then, we can get a naturally induced linear realization of $T'$. But we rather take another linear realization as in Figure \ref{2-1-1} and \ref{2-1-3}, where $\bar{(ii)}$ or $\bar{(iii)}$ means the reverse arcs of (ii) or (iii) (see Figure \ref{inverse}). This is actually another realization of $T'$. Then, $M$-induced link diagram  $L_{T'}$ is isotopic to the $M$-induced link diagram $L_{T}$ or (I)-move connects these two link. Therefore, $L_{T} \in D_{red}$ and $\Sigma(L_{T}) = Y_{T'}=Y_{T}$ by the assumption (see Figure \ref{2-1-1} and \ref{2-1-3}).  
\begin{figure}[h]

\includegraphics[width=12cm,clip]{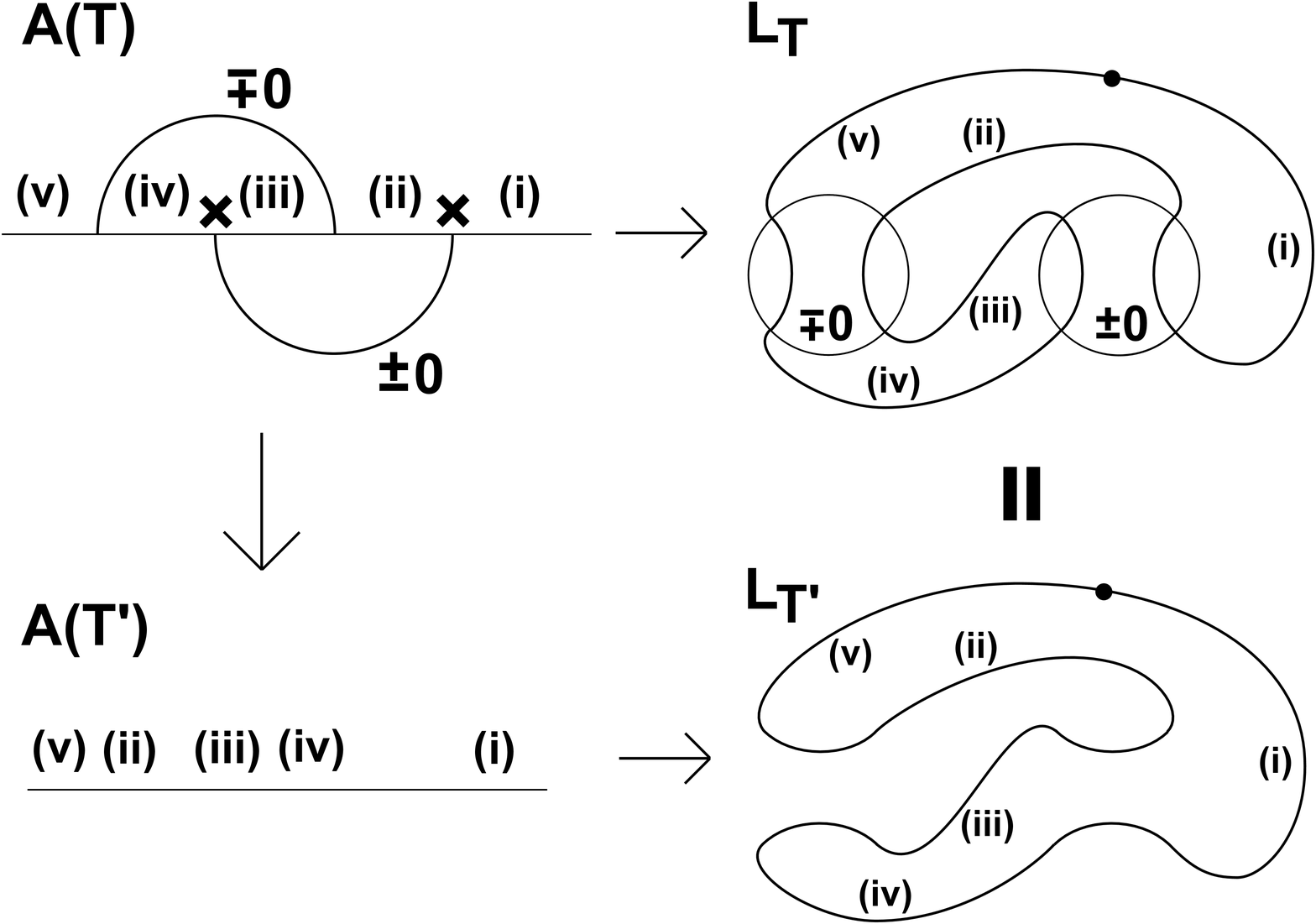} \\
\caption{}\label{2-1-1}
\end{figure}\begin{figure}[h]

\includegraphics[width=12cm,clip]{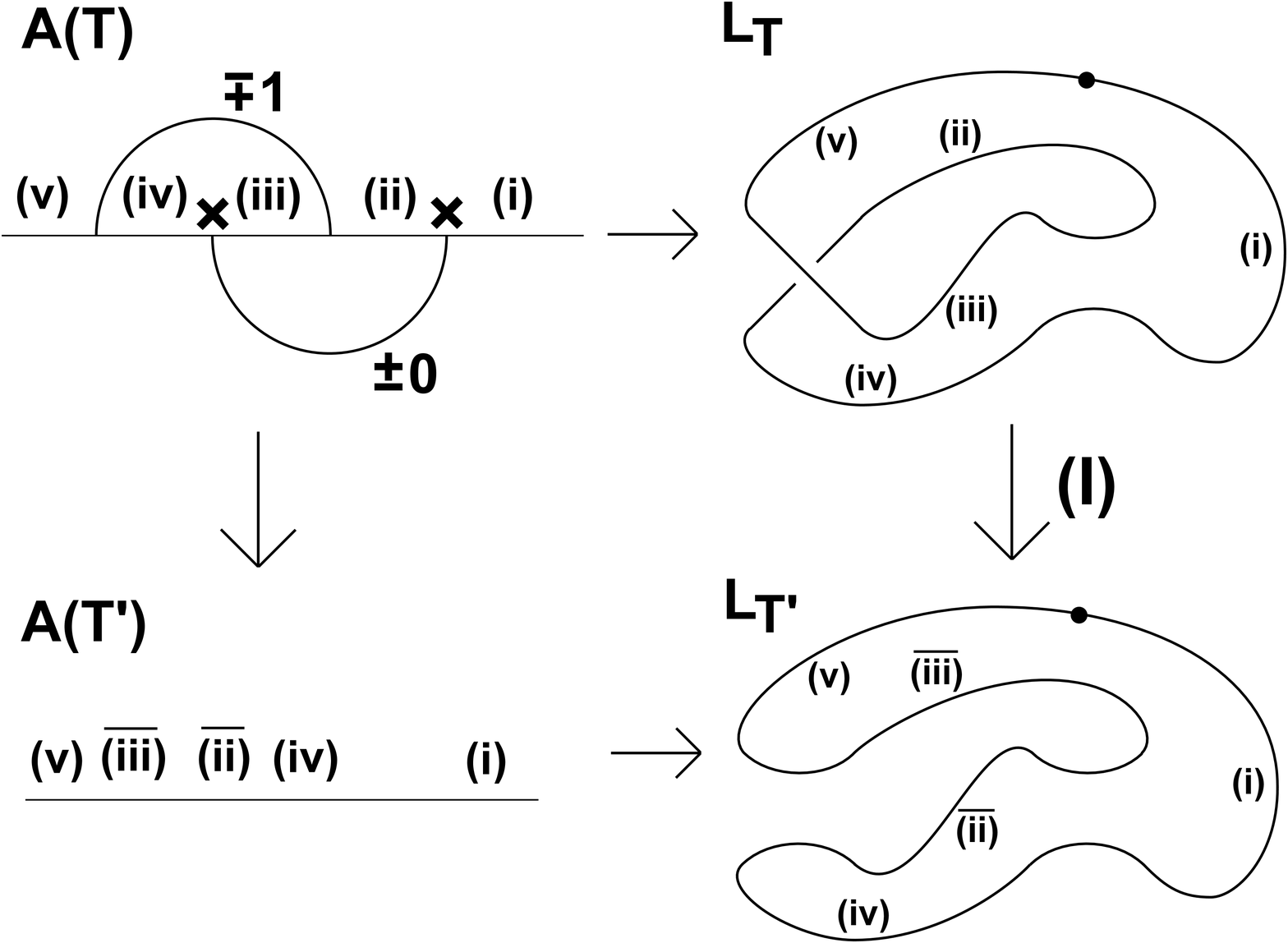} \\
\caption{}\label{2-1-3}
\end{figure}
\begin{figure}[h]

\includegraphics[width=12cm,clip]{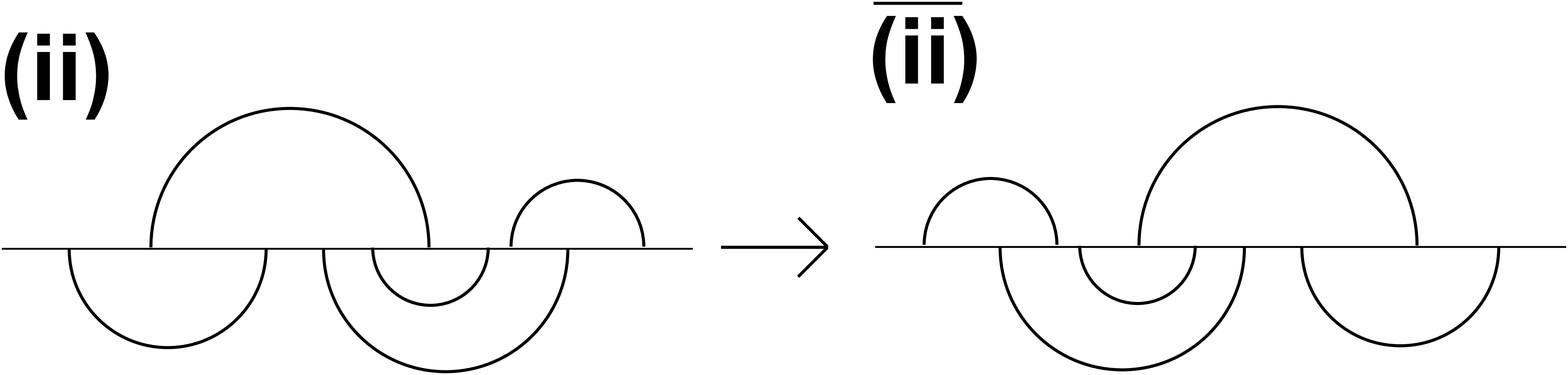} \\
\caption{}\label{inverse}
\end{figure}

\item In this case, for any linear realization $A_{T}$, the natural linear realization $A_{T'}$ is induced. But we should take another linear realization of $T'$ to prove this proposition (see Figure \ref{1-1-1}). Then, $M$-induced link diagram  $L_{T'}$ is isotopic to the $M$-induced link diagram $L_{T}$. Therefore, $L_{T} \in D_{red}$ and $\Sigma(L_{T}) = Y_{T'}=Y_{T}$ by the assumption.

\begin{figure}[h]

\includegraphics[width=12cm,clip]{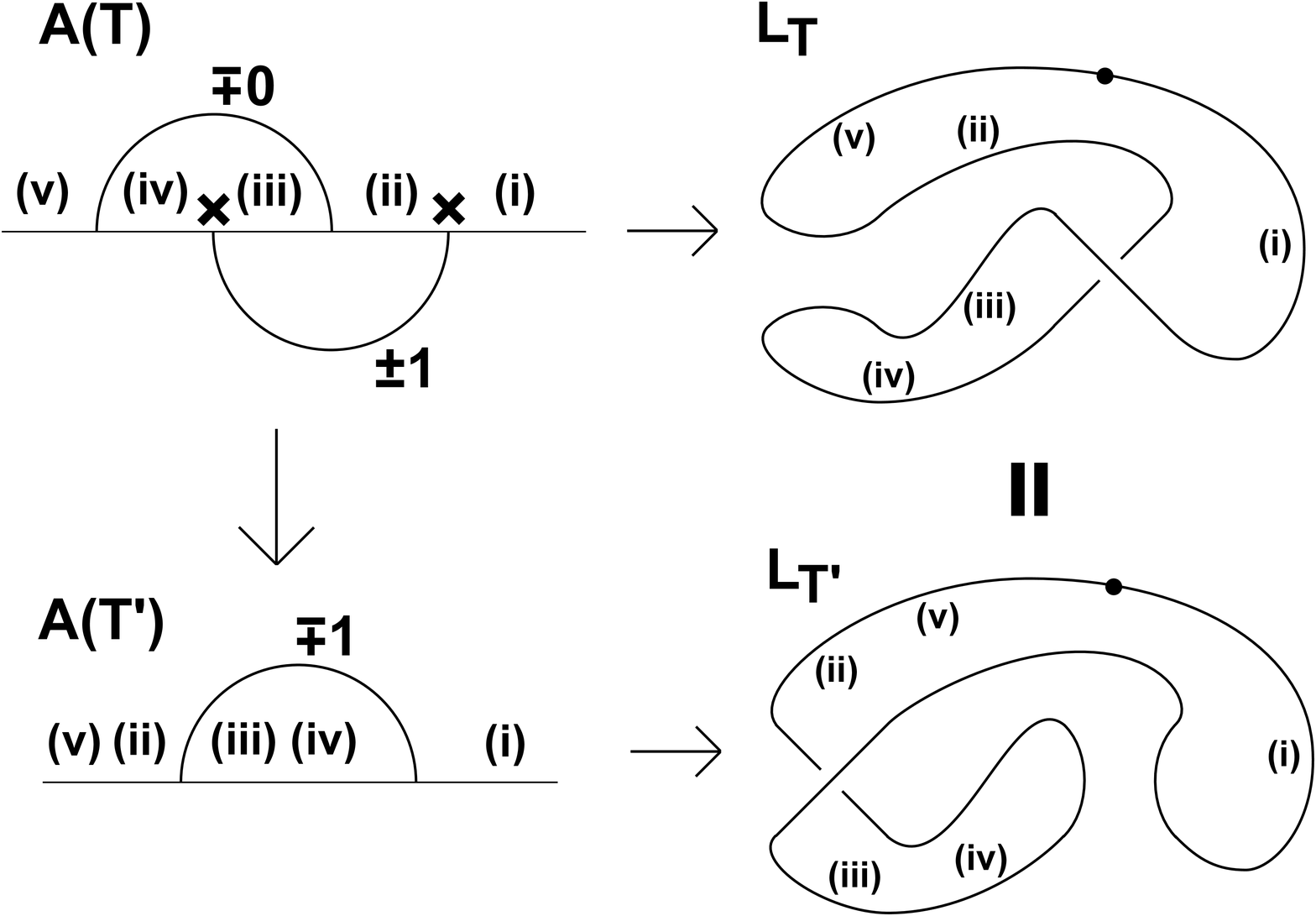} \\
\caption{}\label{1-1-1}
\end{figure}

\item Lastly, we consider the case (4). Any linear realization of $T$ gives the natural linear realization $A_{T'}$ of $T'$ and $M$-induced link diagrams. Then, $L_{T}$ and $L_{T'}$ are connected by one operation (II) (see Figure \ref{3-1-1}). Thus, it holds that $L_{T} \in D_{red}$ and $\Sigma(L_{T}) = Y_{T}$.
\begin{figure}[h]

\includegraphics[width=12cm,clip]{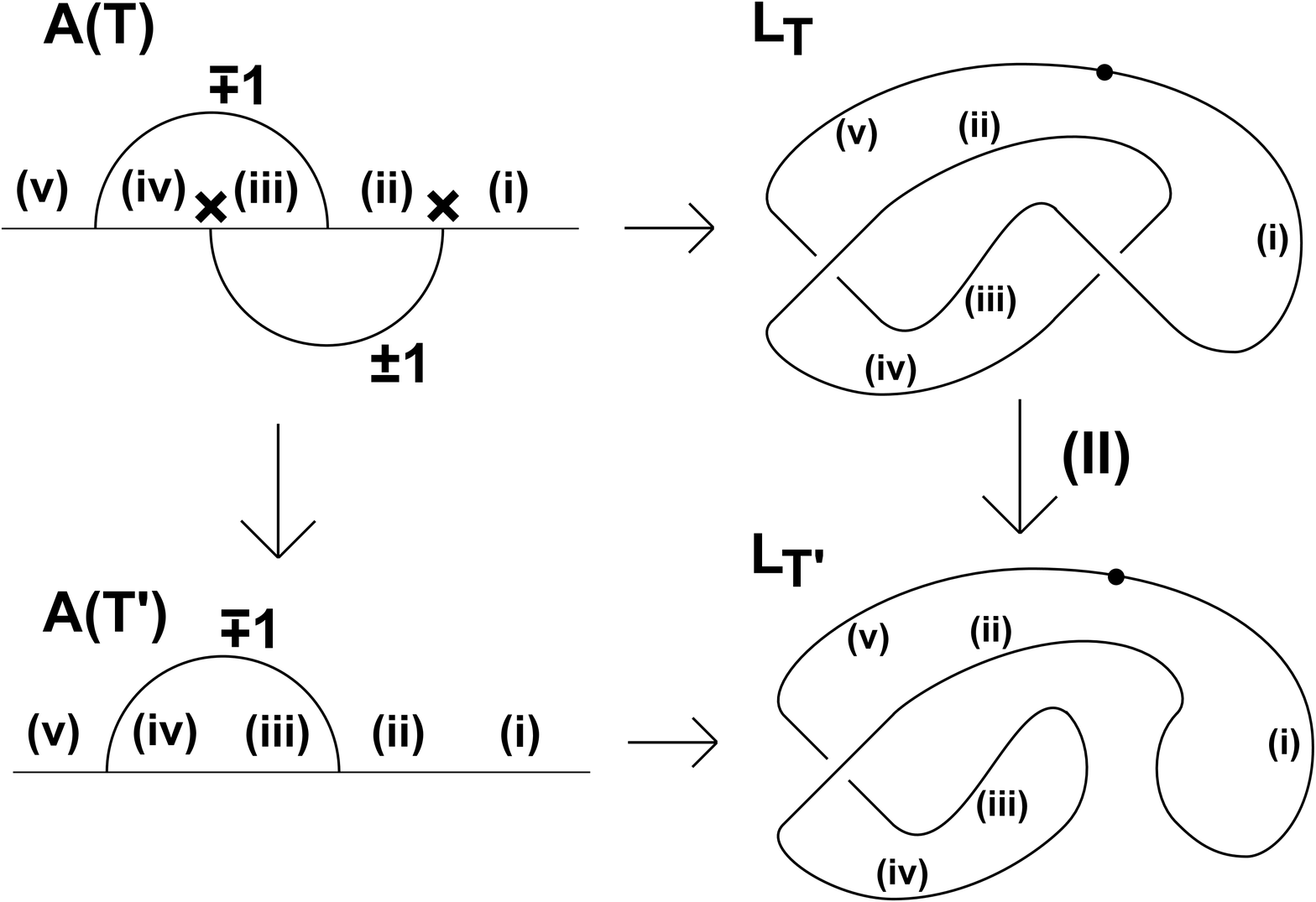} \\
\caption{}\label{3-1-1}
\end{figure}
\end{enumerate}  

$\underline{\mathcal{M}_{\mathcal{T}}\supset \mathcal{M}_{red}}$.
We prove the next claim by induction on the number $|L|$ of the crossing points of the diagram $L$.
\begin{clm}
Let $L \in D_{red}$ and $\Sigma(L) \in \mathcal{M}_{red}$. Then, there exist $\mathbb{T}_{L} \in \mathcal{T}$ and a linear realization $A_{T_{L}}$ of $T_{L}$ such that the $M$-induced link diagram is $L$ and $Y_{T_{L}} = \Sigma(L)$.
\end{clm}
If $|L|=0$, $L$ becomes a disjoint union of unknot diagrams. So we can define $T_{L}$ as finite points with weight $0$.
Next, assume that the proposition holds when $|L| \leqq n$. Take $L \in D_{red}$ with $|L| = n+1$.
Then, we can reduce $L$ by using move (I) or (II) so that a new reduced link $L' \in D_{red}$ has just $n$ crossings. We consider case-by-case.

\begin{enumerate}
\item  In this case, the 1-reducible disk separates $L$ in two parts. We add new arc with weight $\pm 1$ to linear realization $T_{L'}$ and reverse (I) (see Figure \ref{I-1-1}). The $M$-induced link of this new tree $T$ is $L$ and $Y_{T_{L}}=\Sigma(L)$.

\begin{figure}[h]

\includegraphics[width=12cm,clip]{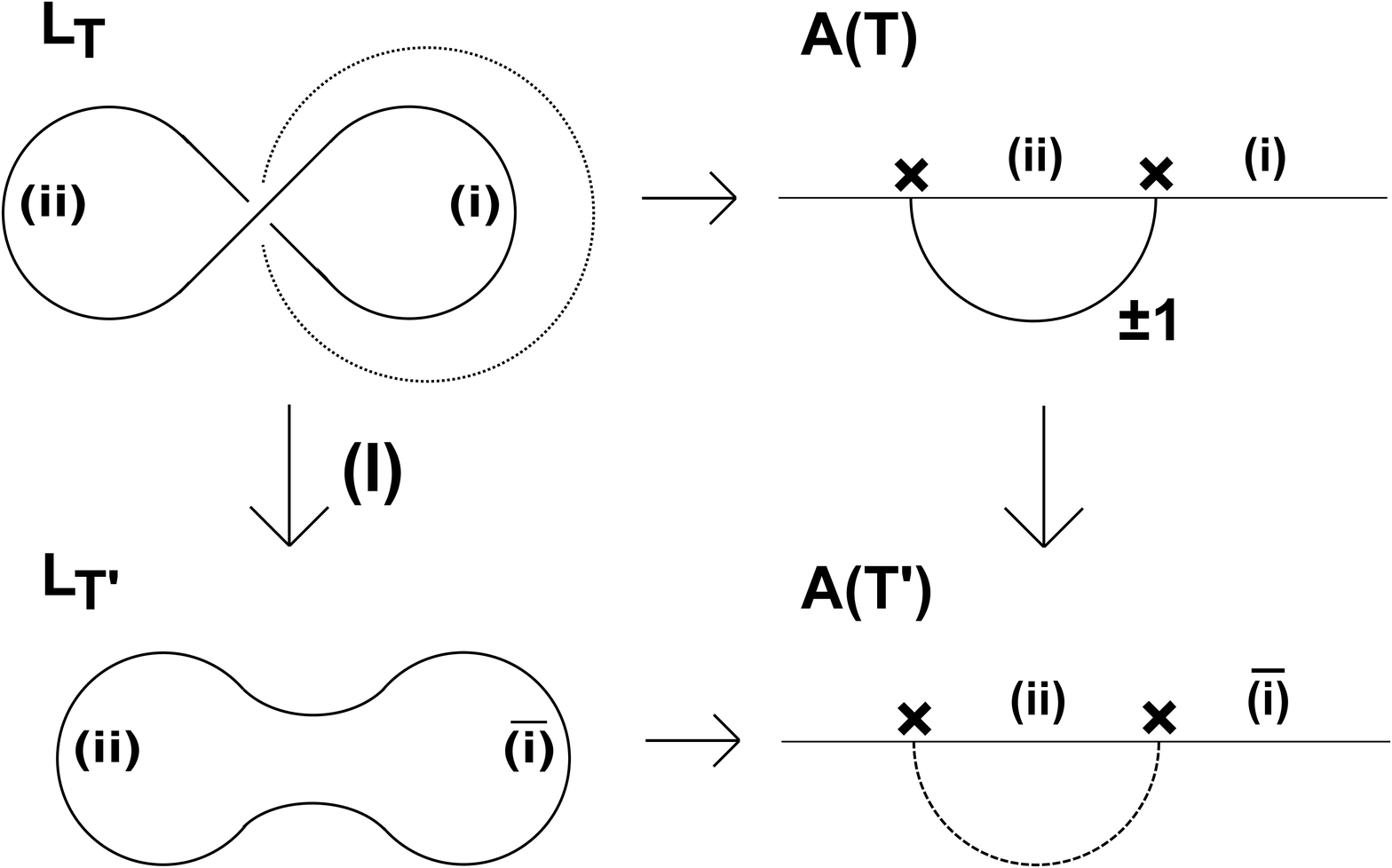} \\
\caption{}\label{I-1-1}
\end{figure}

\item In this case, denote the two crossing $c_{1}$ and $c_{2}$. Assume that the move (II) means to smooth $c_{2}$ as in Figure \ref{II-1-1}. Since we change the new diagram $L'$ into the unknot diagram by smoothing, there are two possible ways to smooth $c_{1}$ (see Figure \ref{II-1-1} (a) and (b)). 

\begin{figure}[h]

\includegraphics[width=10cm,clip]{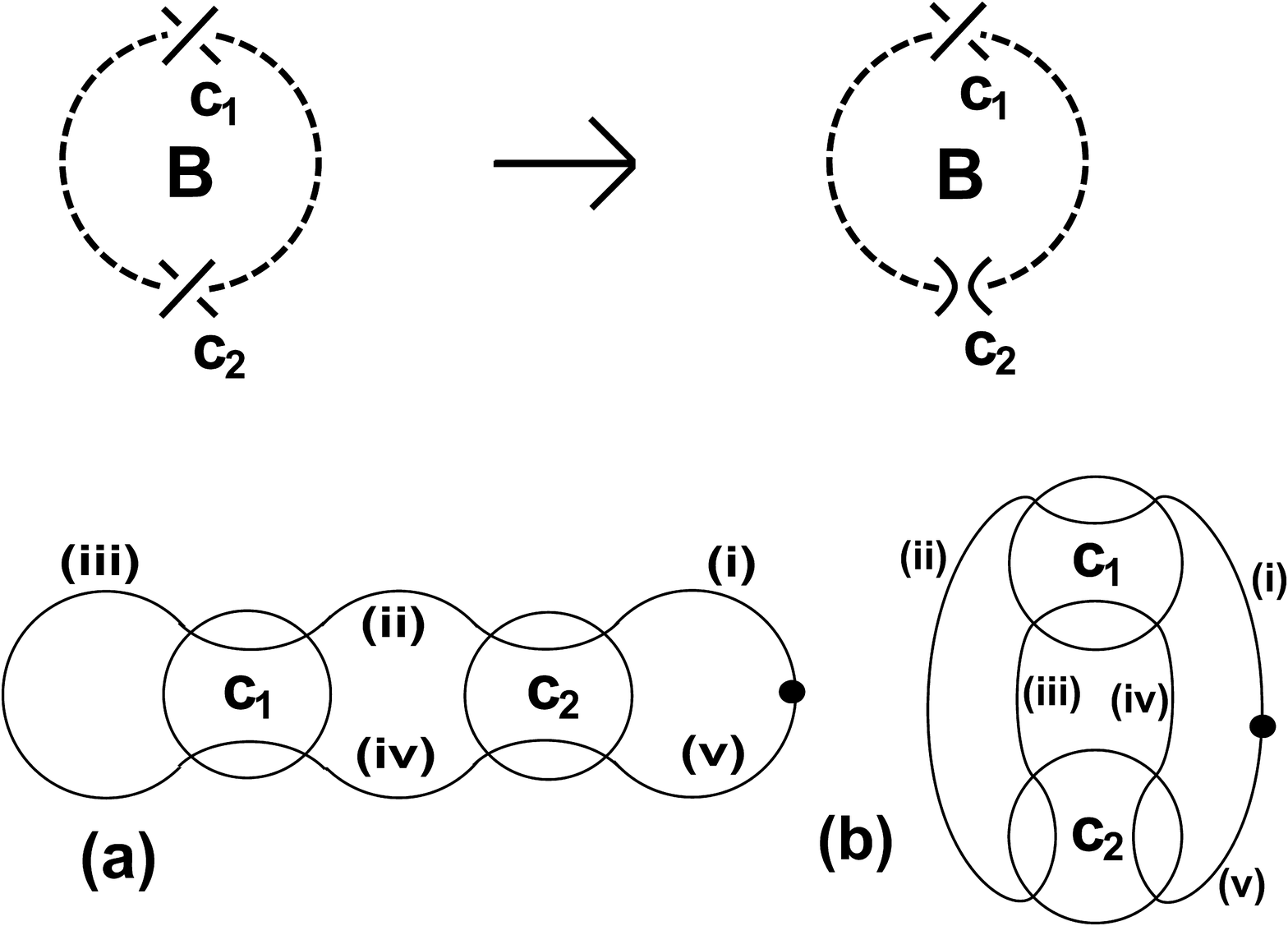} \\
\caption{}\label{II-1-1}
\end{figure}

 \begin{itemize}
 \item In this case, we get a linear realization of $T'$ corresponding $L'$. Since $T'$ is in $\mathcal{T}$, we can define a linear realization of $T \in \mathcal{T}$ as in Figure \ref{II-1-2} if there is at most one arc between (i) and (iii). If there are more than two arcs between (i) and (iii), we should take another linear realization of $T'$ (see Figure \ref{II-1-3}). Then, we can define a linear realization of $T \in \mathcal{T}$ as in Figure \ref{II-1-2} and \ref{II-1-3} and the $M$-induced link diagram is $L$.

\begin{figure}[h]

\includegraphics[width=12cm,clip]{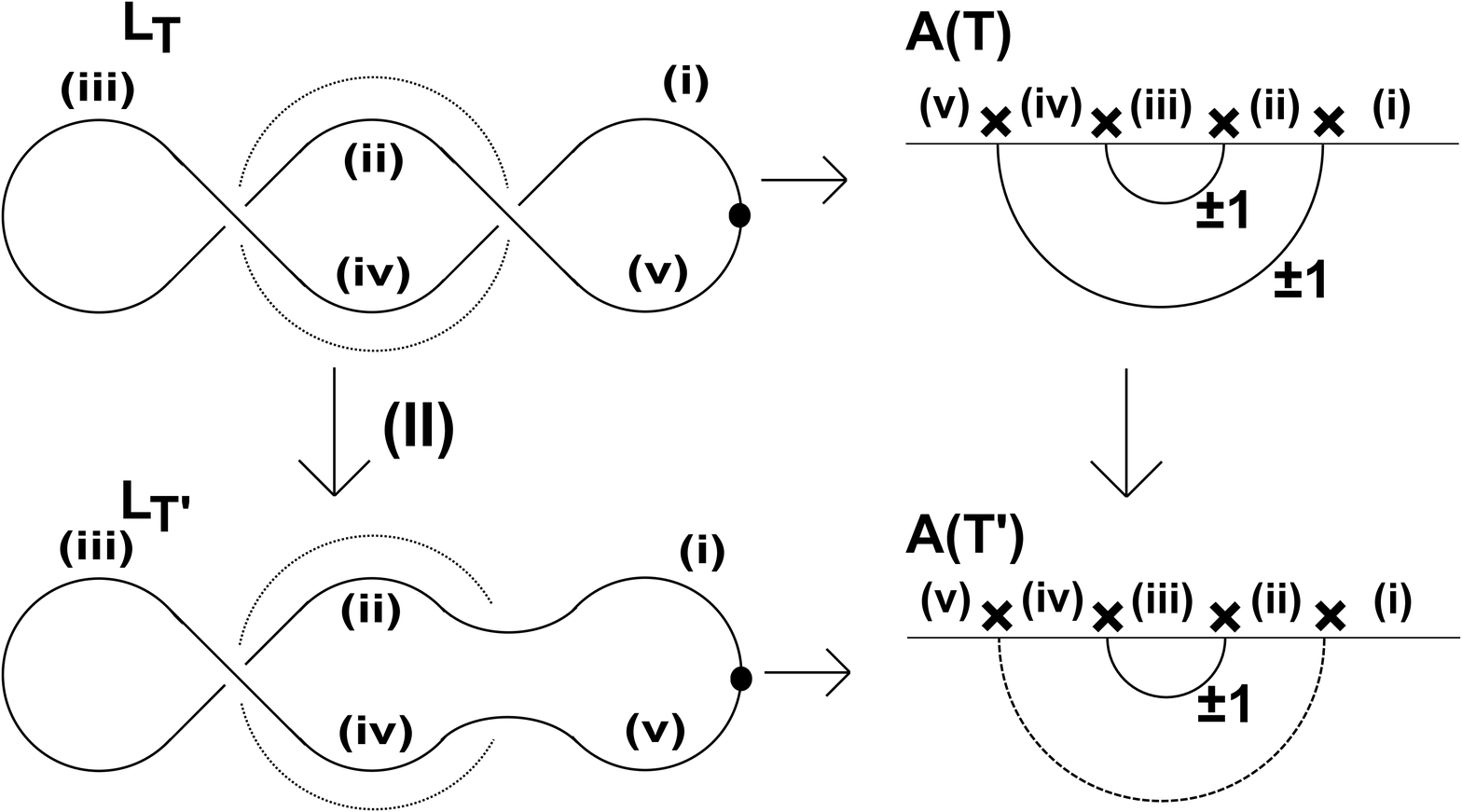} \\
\caption{}\label{II-1-2}
\end{figure}
\begin{figure}[h]

\includegraphics[width=12cm,clip]{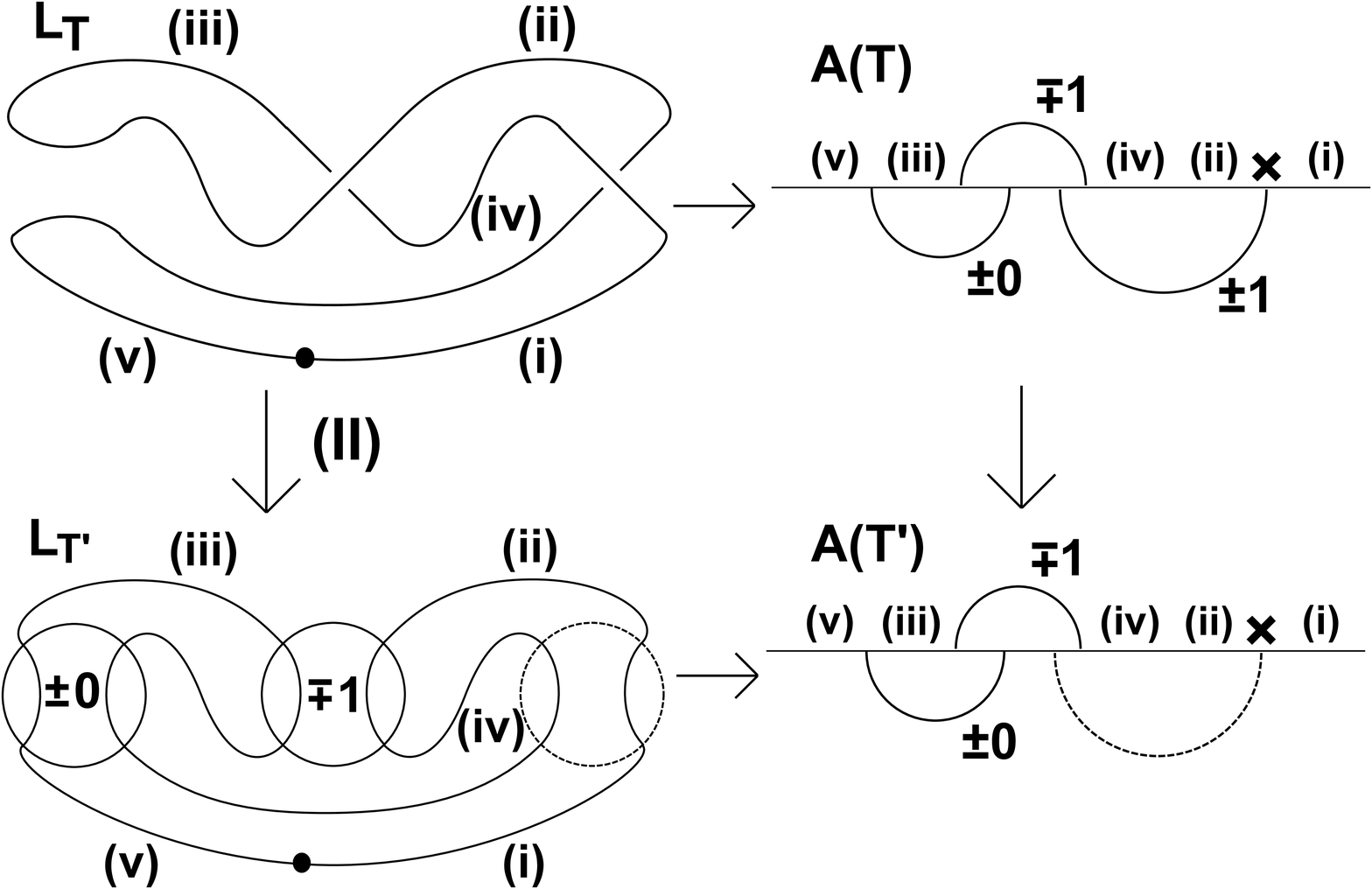} \\
\caption{}\label{II-1-3}
\end{figure}

 \item In this case, we get a linear realization of $T'$ corresponding to $L'$. Since $T'$ is in $\mathcal{T}$, we can define a linear realization of $T$ as in Figure \ref{II-1-4}. Then, $T$ and $T'$ are connected by operation (3). So $T$ is in $\mathcal{T}$ and the $M$-induced link diagram is $L$.
\begin{figure}[h]

\includegraphics[width=12cm,clip]{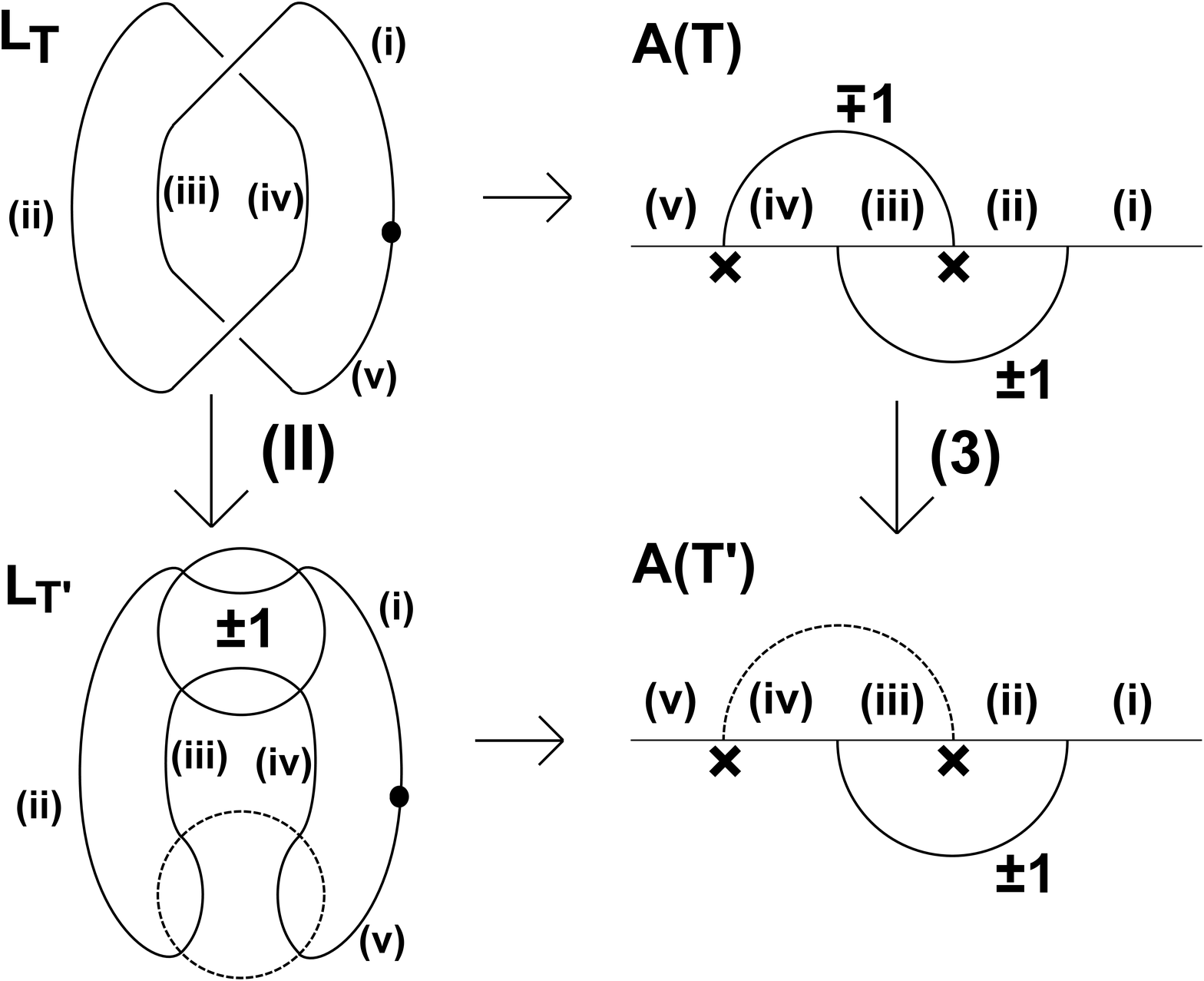} \\
\caption{}\label{II-1-4}
\end{figure}

\end{itemize}  
\end{enumerate}  

\qed
\end{prf}

\section{Proof of Theorem \ref{thmMain}}

\ To prove Theorem \ref{thmMain}, it is enough to prove the next theorem.
We prepare some notations.

\begin{defn}
For a $g \times g$-matrix $A = (a_{ij})$, \textit{expansion signatures of} $\rm{det}(A)$  are signatures of terms $\rm{sgn}(\sigma)a_{1\sigma(1)} \cdots a_{g\sigma(g)}$, where 
$$ \rm{det}(A) = \sum_{\sigma\in S_{g}} \rm{sgn}(\sigma)a_{1\sigma(1)} \cdots a_{g\sigma(g)}.$$
\end{defn}

\begin{defn}
A $g \times g $ matrix $A$ is \textit{effective} if all the non-zero expansion signatures of $\rm{det}(A)$ are constantly positive or constantly negative.
\end{defn}

\begin{thm} \label{thmTvis}
For an alternatingly-weighted tree $T$, if the induced three manifold $Y_{T}$ is a rational homology sphere, then $Y_{T}$ is a strong L-space and a graph manifold (or a connected sum of graphmanifolds).
\end{thm}

\begin{prf}

Let $T$ be an  alternatingly-weighted tree. 
First, we first calculate $|H_{1}(Y;\mathbb{Z})|$. Take an arbitrary ordering on the vertices of $T$. Let $m(v_{i})$ denote the meridian of $K(v_{i})$ for $i = 1, \cdots, g$. Then, these meridians $m(v_{i})$ generate $H_{1}(Y;\mathbb{Z})$ because $L_{T}$ consists of only unknots. All the relations are $\alpha_{i}=0$. This means $|H_{1}(Y;\mathbb{Z})|$ can be calculated by using the following matrix $\rm{Mat}(T)$. Let $w(v_{i})=a(v_{i})/b(v_{i})$ for each vertex $v_{i}$, where $(a(v_{i}),b(v_{i})) = (1,1)$ or $(0,1)$ or $(1,0)$. (Put $1/0 = \infty$.) For each vertex $v_{i}$, the $(i,i)$-components of $\rm{Mat}(T)$ is $\sigma(i)a(i)$, $i = 1,\cdots,g$. For each edge $e$ connecting $i$-th and $j$-th verteces $(i < j)$, the $(i,j)$-th component is $b(i)$ and the $(j,i)$-th component is $b(j)$. The other components are zero. Then, we can calculate $|H_{1}(Y;\mathbb{Z})|$ as the absolute value of the determinant of the matrix $\rm{Mat}(T)$.

Next, take a pointed Heegaard diagram $(\Sigma, \alpha, \beta, z)$ representing $Y_{T}$. Let $L_{T}$ be the induced link from $T$ in $\subset S^3$. Recall that each vertex $v$ of $T$ corresponds to each unknot $K(v)$, and each edge corresponds to linking the two unknots with linking number $\{ \pm 1\}$. So we can take a small arc $c(e)$ for each edge $e$ connecting the two unknots (see Figure \ref{Hdiag}). We can regard the union of the link $L_{T}$ and the arcs $c(e)$ as a spacial graph $G_{T} \subset S^3$. Then, take a small neighborhood of $G_{T}$ and let $\Sigma = \partial G_{T}$. $\Sigma$ is a closed oriented genus $g$ surface, where $g$ is the number of verteces of $T$. (If $T$ is disconnected, we should take tubes connecting these surfaces. This corresponds to connected sums of 3-manifolds.)

Now we assume that $T$ is connected. Then, note that $S^3 \setminus G_{T}$ is a genus $g$ handlebody.  So we can define $\beta$ as its attaching circles. Specifically, each $\beta_{v}$ can be defined near each unknot $K(v)$ as a curve on $\Sigma$ which bounds a disk in $S^3 \setminus G_{T}$ (see Figure \ref{Hdiag}). On the other hand, $\alpha$ curves can be taken as the surgery framings. That is, for each vertex $v$, the weight $\sigma(v)w(v)$ is $\pm 1$ or $\pm 0$ or $\pm \infty$, so $\alpha_{v}$ is defined as a curve on $\Sigma$ with this slope. Take $z$ in $\Sigma \setminus (\alpha \cup \beta)$. Thus, $(\Sigma, \alpha, \beta, z)$ is a pointed Heegaard diagram representing $Y_{T}$. (Note that this diagram is always admissible because $Y$ is a rational homology sphere.) 

\begin{figure}[h]

\includegraphics[width=11cm,clip]{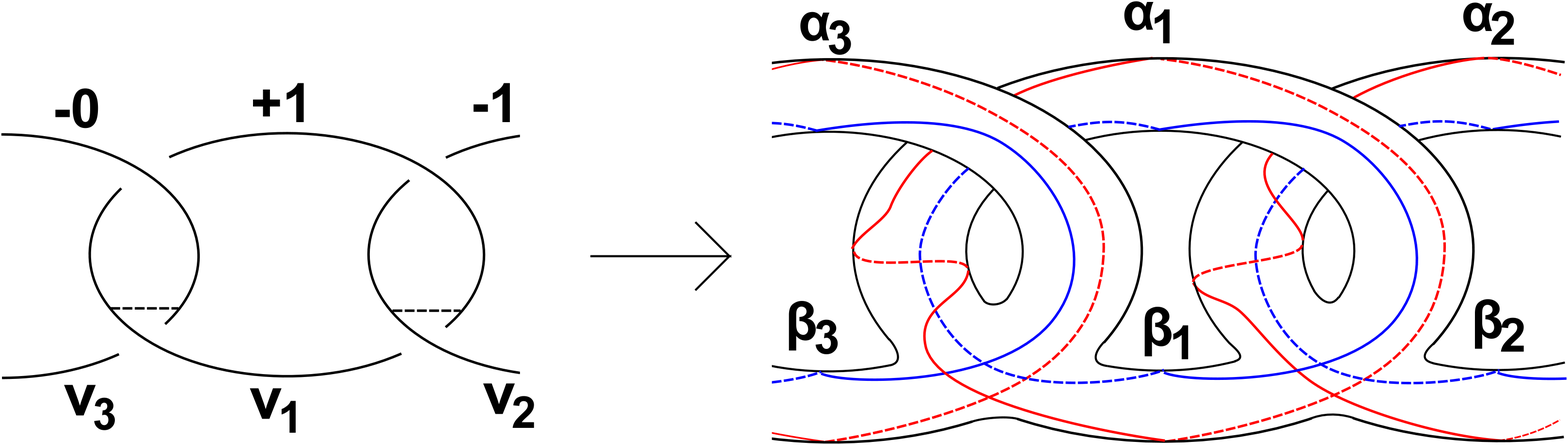} \\
\caption{}\label{Hdiag}
\end{figure}

Now, we compute the number of generators $|\mathbb{T}_{\alpha } \cap \mathbb{T}_{\beta }|$. Note that the number of each local intersection number $\alpha_{i} \cap \beta_{j}$ is the absolute value of the $(i,j)$-component of $\rm{Mat}(T)$. 

It is enough to prove the following proposition.
This proposition proves Theorem \ref{thmMain}. Actually, it implies that:

$$|\mathbb{T}_{\alpha } \cap \mathbb{T}_{\beta }|=|\rm{det}(\rm{Mat}(T))|=|H_{1}(Y;\mathbb{Z})|.$$

Thus, $Y_{T}$ is a strong L-space. 

Finally, we prove the second statement.
To do this, note that by cutting a edge $e$, a connected tree $T$ is decomposed into two trees $T'$ and $T''$. Correspondingly, we can take a torus which decompose $Y_{T}$ into two manifolds with a torus boundary. These manifolds are obviously $Y_{T'}$ and $Y_{T''}$ minus solid tori. By induction of the number of the vertex of $|T|$, we finish the proof.
\qed
\end{prf}

\begin{prop} \label{proplast}
For an alternatingly-weighted tree $T$, $\rm{Mat}(T)$ is an effective matrix.
\end{prop}
\begin{prf4}

we prove this proposition by induction on the number $g$ of vertices of $T$.
If $g=1$, it is trivial. If $g=2$, it is easy because $T$ has alternating weight.
We assume that $g-1$ and $g-2$ cases are proved.

First, fix one univalent vertex $v_{1}$ and denote the next vertex $v_{2}$. Let $T'$ denote the tree without the vertx $v_{1}$ and the unique edge connecting $v_{1}$. Similarly let $T''$ denote the tree without $v_{1}$ and $v_{2}$ and the edges connecting $v_{1}$ and connecting $v_{2}$. Then, we get two another matrices $\rm{Mat}(T')$ and $\rm{Mat}(T'')$. By the above assumptionm, $\rm{Mat}(T')$ and $\rm{Mat}(T'')$ have constant expansion signatures. Denote them $\rm{sgn}(T')$ and $\rm{sgn}(T'')$.
Moreover, these signature satisfies $\rm{sgn}(T') = \sigma(2)\rm{sgn}(T'')$ because $T'$ has also an alternating weight.

We put $w(1)=a(1)/b(1)$ and $w(2)=a(2)/b(2)$. Note that $\sigma(1)=-\sigma(2)$. So $\rm{det}(\rm{Mat}(T))$ satisfies the following equation.

$$\rm{det}(\rm{Mat}(T))=\sigma(1)a(1)\rm{det}(\rm{Mat}(T'))-b(1)b(2)\rm{det}(\rm{Mat}(T'')).$$

Then, the expansion signatures are constant because $$\sigma(1)\rm{sgn}(T')=-\sigma(2)^2\rm{sgn}(T'')=-\rm{sgn}(T'').$$
\qed
\end{prf4}

\begin{prf5}

Theorem \ref{thmRedBrm}, Theorem \ref{thmTred} and Theorem \ref{thmTvis} imply.
\qed
\end{prf5}
%

%
\section*{acknowledgement}
I would like to express my deepest gratitude to Prof. Kohno who provided helpful comments and suggestions. I would also like to express my gratitude to my family for their moral support and warm encouragements.

\end{document}